# Free-fermionic ice, the Yang-Baxter equation and skein relations

Chavdar Lalov

Stanford University, CA

**Abstract**

In his famous ASM paper, Kuperberg uses a skein relation to give an algebraic proof of a Yang-Baxter equation where the Boltzmann weights satisfy the field-free condition. In this paper, we use Kuperberg's techniques to give proofs of a few Yang-Baxter equations where the Boltzmann weights satisfy the free-fermionic condition. In particular, we use skein relations to prove the Yang-Baxter equation for Gamma-Gamma ice which is a free-fermionic six-vertex model introduced by Brubaker, Bump and Friedberg.

## Contents

# 1 Introduction

In his famous paper [K] Kuperberg uses the field-free six-vertex model to give a short second proof of the Alternating Sign Matrix Conjecture (after the purely combinatorial one by Zeilberger [Z]). He uses a bijection between alternating sign matrices and states of the six-vertex model to bring the problem into the realm of solvable lattice models. Building on work by Baxter [B], and Korepin and Izergin [I2, BIK], Kuperberg computes the partition function of the respective six-vertex model, thus finding the number of alternating sign matrices.

A key ingredient in Kuperberg's proof is the field-free Yang-Baxter equation. It is rather standard to use a computer when proving the Yang-Baxter equation as multiple cases must be checked. Kuperberg, however, takes a different approach. He treats the Yang-Baxter equation as a knot and finds a skein relation that decomposes the Boltzmann weights. This allows him to give a quick algebraic proof of the Yang-Baxter equation by relating it to the Temperley-Lieb algebra and the quantum group $U_q(sl(2))$. What is intriguing here is not so much that Kuperberg avoided the use of a brute force computer calculation, but rather that he found a surprising connection with knot theory which brings a new perspective to the Yang-Baxter equation. One may argue that the method Kuperberg used to prove the field-free Yang-Baxter equation is just as interesting and insightful as the result itself.

This paper is inspired by the following question: Can Kuperberg's method for proving the field-free Yang-Baxter equation be used to prove Yang-Baxter equations where the Boltzmann weights do not satisfy the field-free condition. We present several results that give an affirmative answer to this question. It is worth noting that although some of our techniques resemble those of Kuperberg, others are very different. Our main insight is that different skein relations can be used to prove important free-fermionic Yang-Baxter equations. One reason for why free-fermionic six-vertex models and their respective Yang-Baxter equations are important is their application in proving Tokuyama formulas for different reductive groups [BBF], [I1].

We strongly feel that the methods we developed to prove Yang-Baxter equations in this paper are just as interesting as the results themselves. Indeed, some of the theorems we present already have computer-assisted proofs in other articles, however, those proofs do not reveal all of the hidden mathematical structure. We suspect that the combinatorial techniques and skein relations we discovered might have further applications in knot theory, e.g. knot invariants, as well as in solvable lattice models although this direction is yet to be studied.

In Section 2 of this paper we present a proof of a certain free-fermionic Yang Baxter equation and the respective Boltzmann weights are given in Table 1. The idea we used is rather surprising. Firstly, we prove two different field-free Yang-Baxter equations by using techniques similar to Kuperberg's. The Boltzmann weights for these Yang-Baxter equations are given in Tables 2 and 3. Then we use a combinatorial lemma to mix the two field-free Yang-Baxter equations into a free-fermionic one. The lemma states that the $a_1$ Boltzmann weight and the $a_2$ Boltzmann weight cannot appear simultaneously in any instance of the Yang-Baxter equation. It is worth noting that the cap weights used in the skein relations in this section have very nice properties such as the ability to straighten curved lines. This is a standard property in knot theory.

In Section 3 we further develop our techniques from Section 2 to prove another free-fermionic Yang-Baxter equation from a paper by Brubaker, Bump and Friedberg [BBF]. The Boltzmann weights in Section 3 are known as Gamma-Gamma ice [BBF,I1] and can be found in Table 3. The section contains another important lemma that allows us to move constants between the $b_1$ and $b_2$ Boltzmann weights.

Section 4 is optional and can be skipped as it is not needed to read the other sections. In it we explain the connection between the field-free Yang-Baxter equation of Kuperberg and the two field-free Yang-Baxter equations in Section 2.

In Sections 5, 6 and 7 we prove a more general form of the Yang-Baxter equation for Gamma-Gamma ice which again can be found in the article by Brubaker, Bump and Friedberg [BBF]. The Boltzmann weights can be found in Table 8. It is interesting that our methods from Sections 2 and 3 do not seem to work in this case. As a result, we present a completely new method that differs quite a bit from Kuperberg's. The method depends on an intriguing connection between the fish equation and the Yang-Baxter equation. The skein relations we find are very different from the ones in Sections 2 and 3. It is important to note that Section 7 gives a second proof of the result in Section 3.

## 2 A free-fermionic Yang-Baxter equation

In this section we show that one can use Kuperberg's methods in [K] to prove a Yang-Baxter equation where the Boltzmann weights satisfy the free-fermionic condition.

The Bolztmann weights we use are

| $a_1$ | $a_2$ | $b_1$ | $b_2$ | $c_1$ | $c_2$ |
|---|---|---|---|---|---|
| 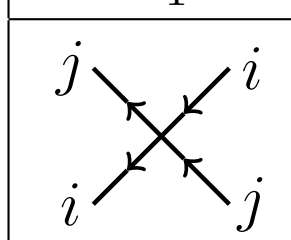 | 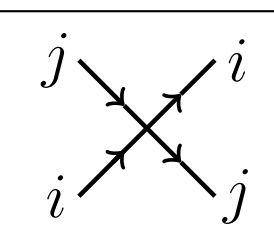 | 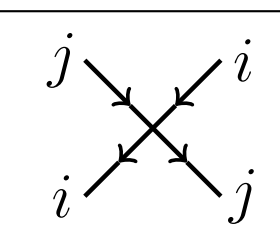 | 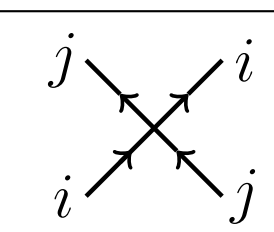 | 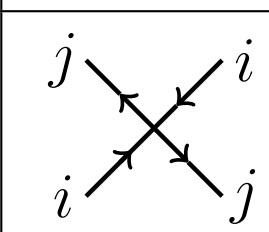 | 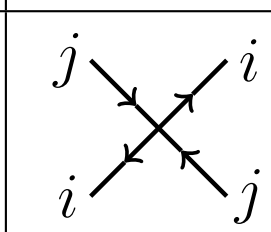 |
| $tz_i + z_j$ | $tz_j + z_i$ | $\sqrt{-t}(z_i - z_j)$ | $\sqrt{-t}(z_i - z_j)$ | $(t+1)z_i$ | $(t+1)z_j$ |

Table 1

We will refer to the edge orientations of the vertices as either left or right. For example, $a_1$ consists of 4 left arrows while $b_1$ consists of 2 left and 2 right arrows. Each Boltzmann weight is parameterised by the rows $i$ and $j$. The parameters $z_i$ and $z_j$ are called *spectral parameters* while the parameter $t$ is called a *deformation parameter*. We call an edge *incoming* if it points towards the vertex and call it *outgoing* otherwise. Notice that each vertex in Table 1 has two incoming and two outgoing arrows. The weights satisfy the free-fermionic condition:

$$a_1 a_2 + b_1 b_2 = c_1 c_2.$$

In this section we prove the weights in Table 1 satisfy the Yang-Baxter equation.

**Theorem 2.1.** *For any way the external edges are fixed, the Boltzmann weights in Table 1*

*satisfy*

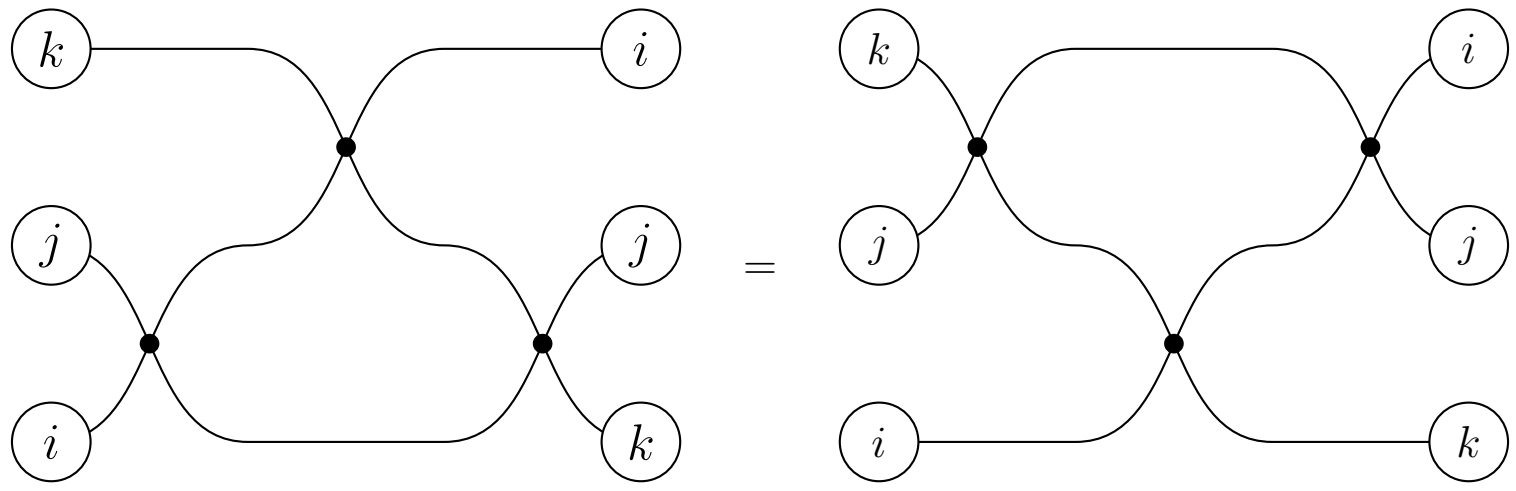

*where on both sides of the equation we sum over all possible configurations of the internal edges.*

**Example 2.1** (Yang-Baxter equation)**.** *We present an instance of the Yang-Baxter equation for one choice of the external edges. Recall that it is mandatory for the external edges to have the same orientation on the two sides of the equation:*

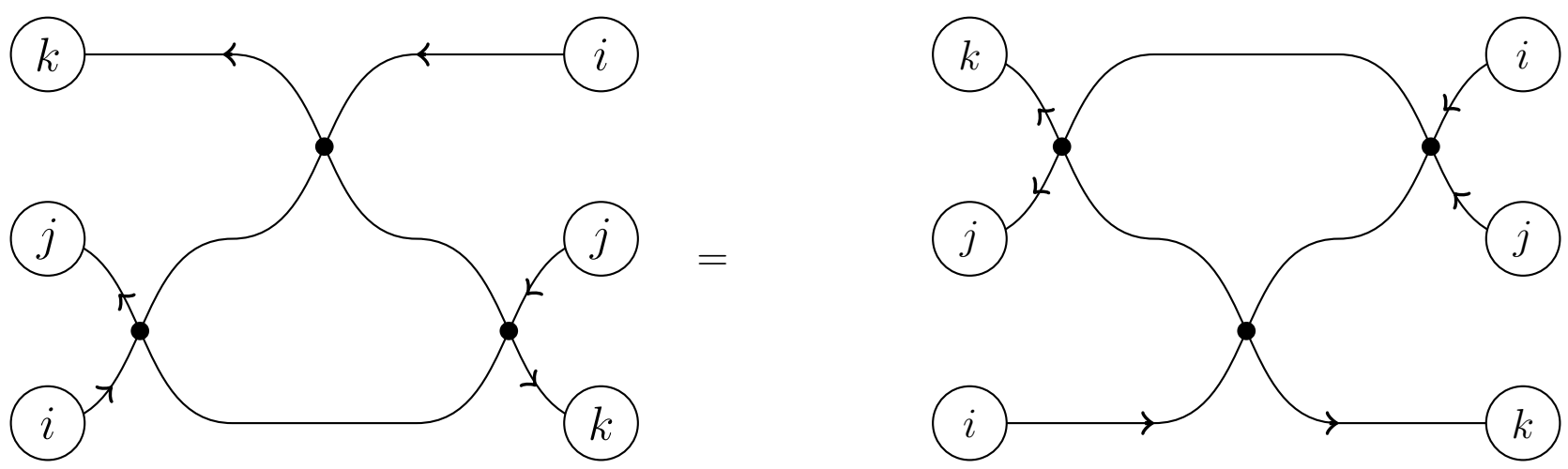

*Next we fill in the interior in all possible ways and sum over the weights of all configurations on both sides of the equation. Note that the weight of each configuration is just the product of the Boltzmann weights in it.*

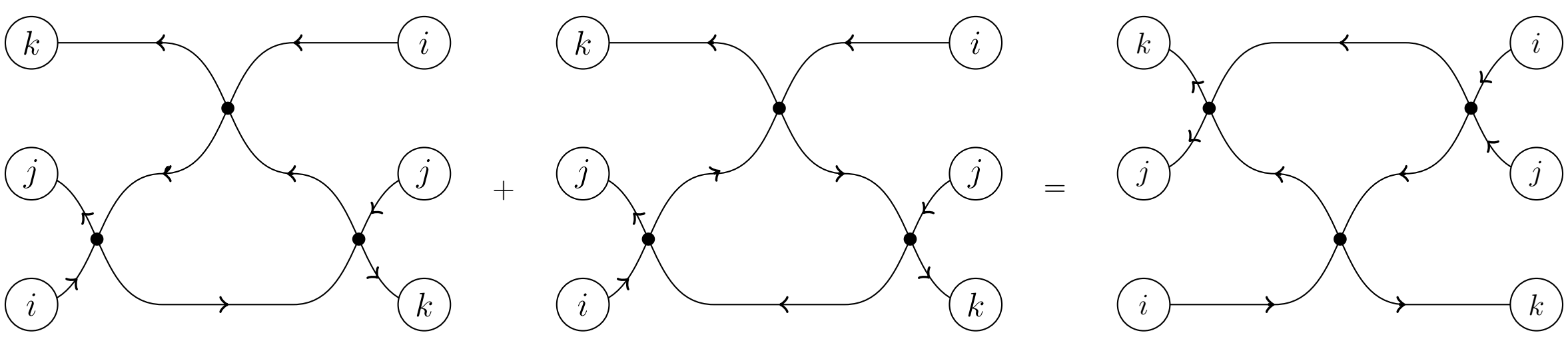

*The identiy is equiavlent to*

$$(tz_i + z_k)(t+1)z_i(t+1)z_j + (t+1)z_i\sqrt{-t}(z_i - z_j)\sqrt{-t}(z_j - z_k) = (tz_j + z_k)(tz_i + z_j)(t+1)z_i$$

*which is not difficult to check.*

Note that all non-trivial instances of the Yang-Baxter equation have three incoming external edges and three outgoing external edges. Otherwise both sides of the Yang-Baxter equation will equal zero. This is because each vertex in Table 1 has two incoming and two outgoing edges. As a result, each non-trivial instance of the Yang-Baxter equation has an even number of external edges pointing to the right and an even number of external edges pointing to the left.

To prove Theorem 2.1 we introduce an important lemma.

**Lemma 2.2.** *For any non-trivial instance of the Yang-Baxter equation in which there are at least four external edges pointing to the left, there is no $a_2$ state. Similarly, for any non-trivial instance of the Yang-Baxter equation in which there are at least four external edges pointing to the right, there is no $a_1$ state. In particular, $a_1$ and $a_2$ cannot appear simultaneously in any non-trivial instance of the Yang-Baxter equation.*

*Proof.* We prove the first statement. Assume for the sake of contradiction that there is an instance of the Yang-Baxter equation with at least four external edge pointing to the left where there is an $a_2$ state present in the equation. Notice that we can either have four or six external edges pointing to the left. We cannot have five as that would imply we have at least four incoming or four outgoing external edges. The case with six external edges pointing to the left is impossible as $a_2$ has only edges pointing to the right. Thus, we must have exactly four external edges pointing to the left and exactly two pointing to the right. The two external edges pointing to the right must be part of the $a_2$ state. Moreover, one of the external edges pointing to the right must be incoming and one must be outgoing. Without loss of generality assume that the $a_2$ state is present on the right-hand side of the Yang-Baxter equation. Then the above discussion forces the following configuration

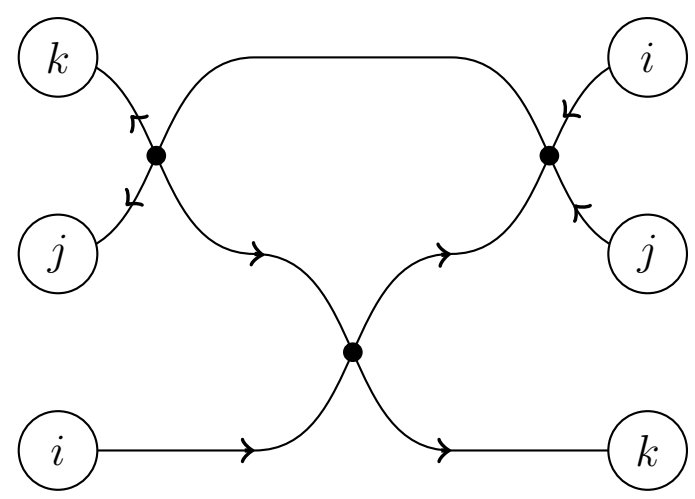

Figure 2: *Forced configuration assuming presence of an $a_2$ state and exactly four external edges pointing to the left.*

The configuration in Figure 2 is not valid as there are vertices with three incoming and with three outgoing edges. We have reached a contradiction.

The second statement of the lemma follows from the first. This is because if there exists a configuration with at least four external edges pointing to the right and an $a_1$ state, then by rotating the picture by 180 degrees we would get a valid configuration with at least four external edges pointing to the left and an $a_2$ state. The latter leads to a contradiction.

Finally, we cannot have an $a_1$ and an $a_2$ state simultaneously in any non-trivial instance of the Yang-Baxter equation because we will always have either four external edges pointing to the right or four external edges pointing to the left. □

Thus, Lemma 2.2 implies that the states $a_1$ and $a_2$ divide the instances of Yang-Baxter equation into two groups. One group never has an $a_2$ state while the other never has an $a_1$ state.

Next we use Kuperberg's techniques from [K] to prove that both rows of Boltzmann weights in Table 2 satisfy the Yang-Baxter equation. The weights in Table 2 differ from the weights in Table 1 by only one weight (either $a_1$ or $a_2$).

| | $a_1$ | $a_2$ | $b_1$ | $b_2$ | $c_1$ | $c_2$ |
|---|---|---|---|---|---|---|
| | 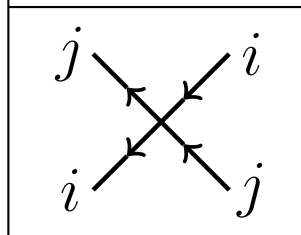 | 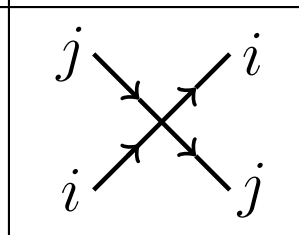 | 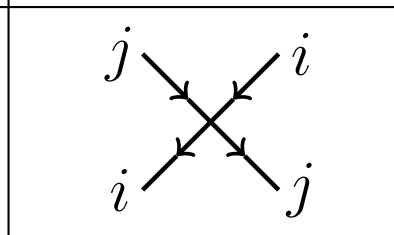 | 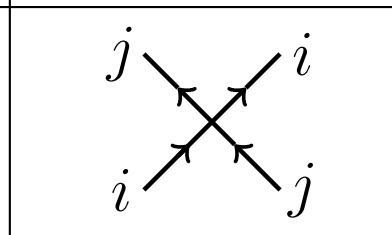 | 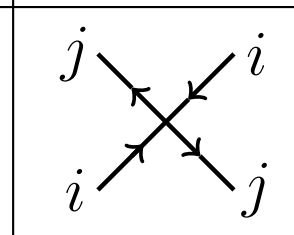 | 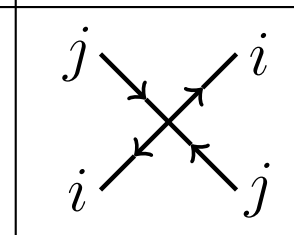 |
| Row 1 | $tz_i + z_j$ | $tz_i + z_j$ | $\sqrt{-t}(z_i - z_j)$ | $\sqrt{-t}(z_i - z_j)$ | $(t+1)z_i$ | $(t+1)z_j$ |
| Row 2 | $tz_j + z_i$ | $tz_j + z_i$ | $\sqrt{-t}(z_i - z_j)$ | $\sqrt{-t}(z_i - z_j)$ | $(t+1)z_i$ | $(t+1)z_j$ |

Table 2: *Row 1 an Row 2 are two different combinations of Boltzmann weights that differ only by the value of $a_1$ or $a_2$. Both rows have field-free Boltzmann weights.*

**Proposition 2.3.** *For any way the external edges are fixed, the Boltzmann weights in both Row 1 and Row 2 of Table 2 satisfy:*

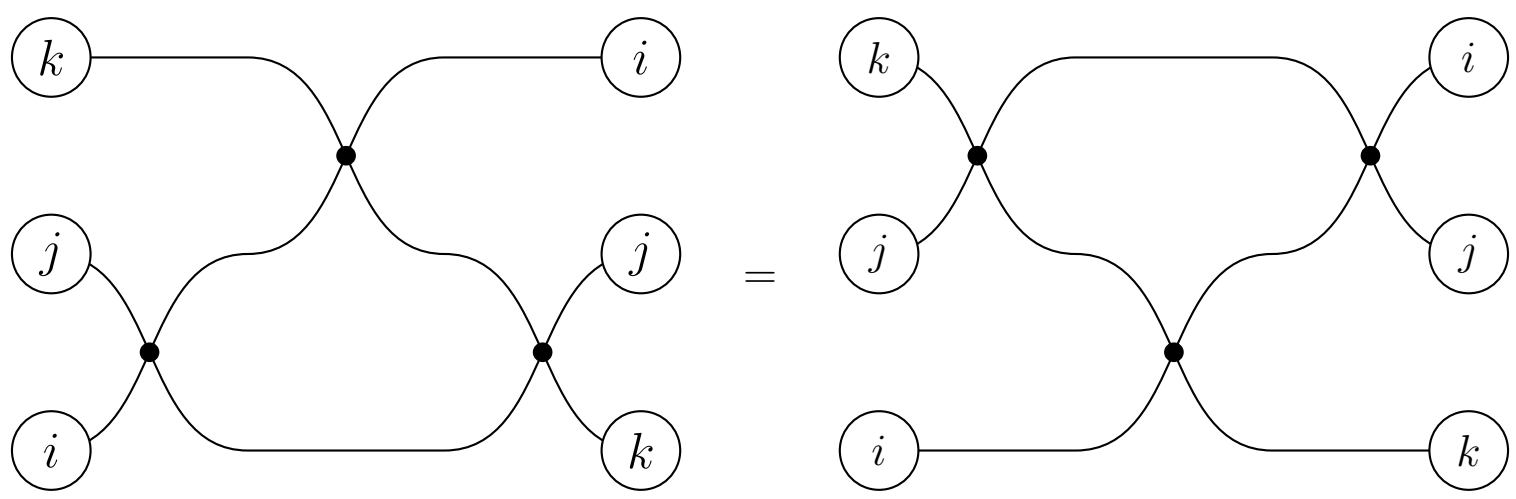

*Proof of Proposition 2.2.* We divide the proof into two cases, one for each row of Boltzmann weights in Table 2.

**First case:** We begin by proving the Boltzmann weights in Row 1 of Table 2 satisfy the Yang-Baxter equation.

To do this we use a skein relation that decomposes the Boltzmann weights.

Define the following cap weights:

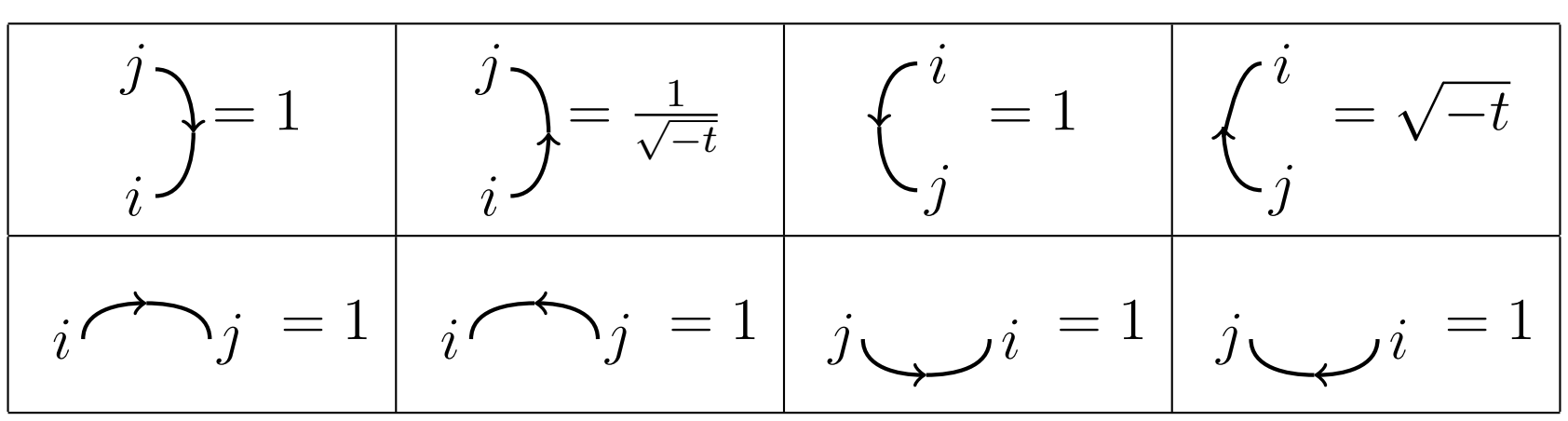

**Lemma 2.4.** *The Boltzmann weight from Row 1 in Table 2 satisfy the following skein relation*

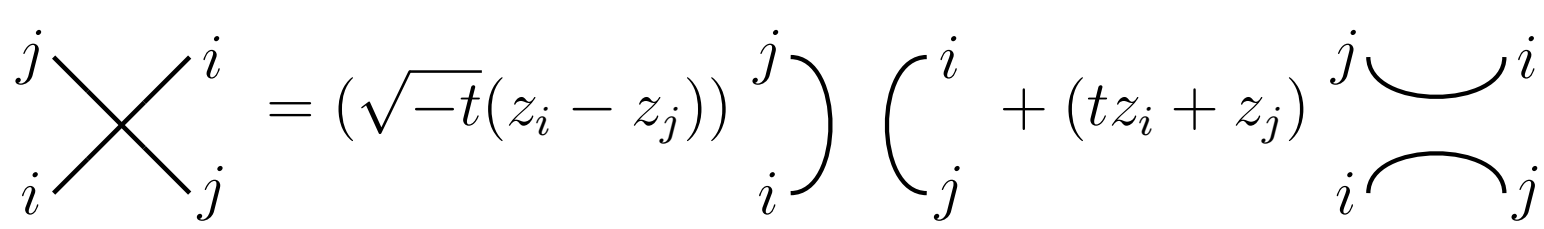

(1)

*Proof.* The skein relation is verified by direct computation. □

The cap weights, we defined, have a few more important properties which are straightforward to check:

**Lemma 2.5.** *The following identities hold:*

***Property 1:***

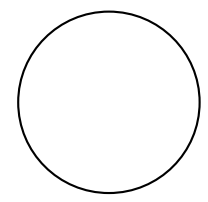

$$= \sqrt{-t} + \frac{1}{\sqrt{-t}} \tag{2}$$

*where we are summing over the two possible orientations of the loop.*

***Property 2:***

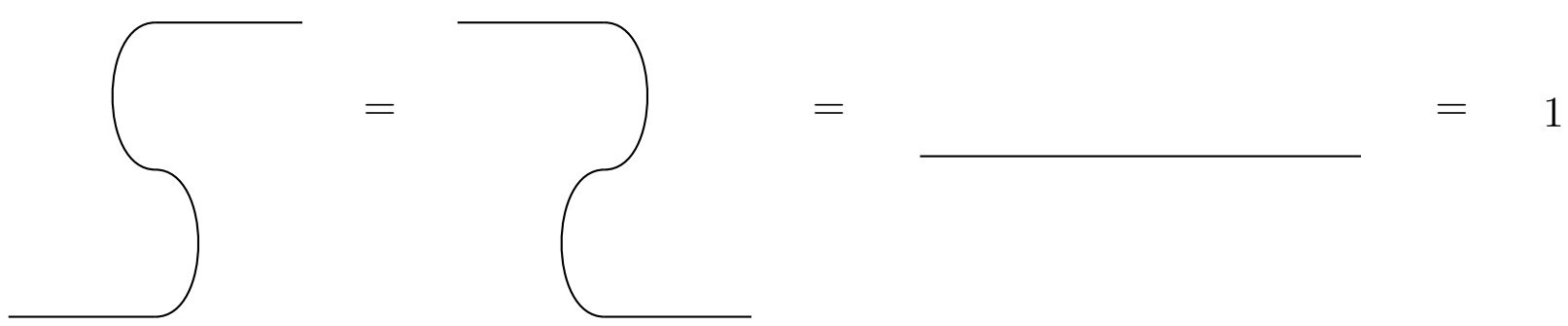

(3)

*where the equation holds for both orientation of the curvy line and the straight line.*

***Property 3:***

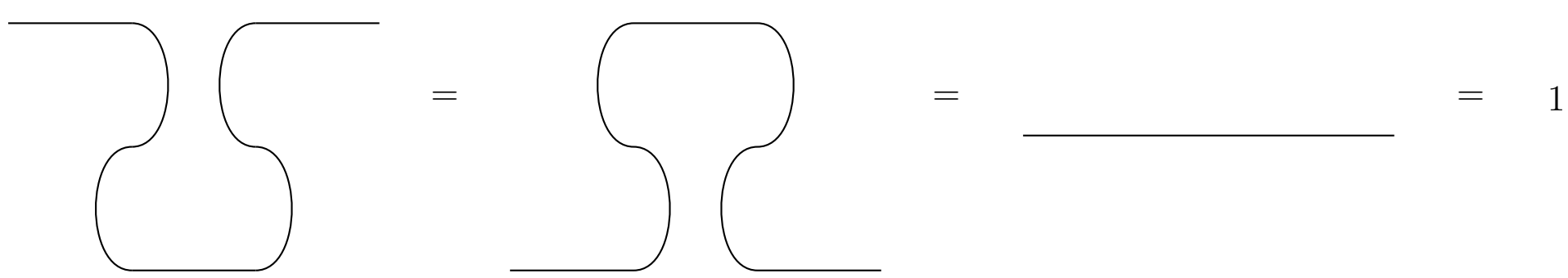

(4)

*where the equation holds for both orientations of the curvy lines and the straight line.*

Note that the second and third property of Lemma 2.5 imply that we can straighten curvy lines which is typical in applications of the six-vertex model to knot theory. This is also an important property of the graphic representation of the Temperley-Lieb algebra.

By using the skein relation (1) we can expand the left-hand side of the Yang-Baxter equation into eight terms.

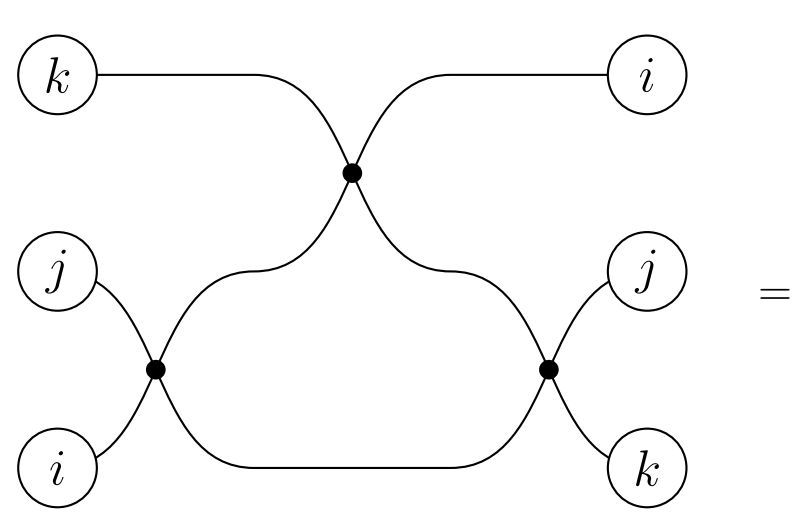

$= \quad (tz_i + z_j)(tz_i + z_k)(tz_j + z_k)$ 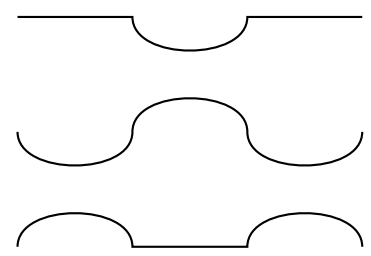 $+ \quad (tz_i + z_j)\sqrt{-t}(z_i - z_k)\sqrt{-t}(z_j - z_k)$ 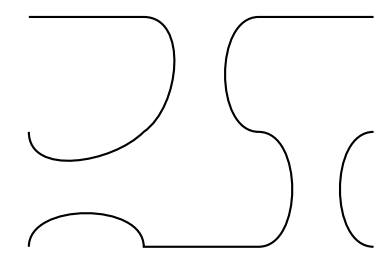

$+ \quad \sqrt{-t}(z_i - z_j)\sqrt{-t}(z_i - z_k)(tz_j + z_k)$ 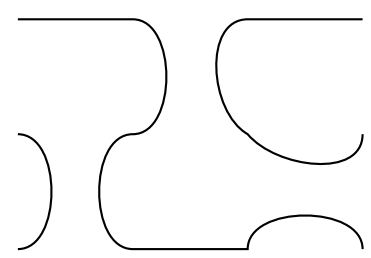 $+ \quad (tz_i + z_j)\sqrt{-t}(z_i - z_k)(tz_j + z_k)$ 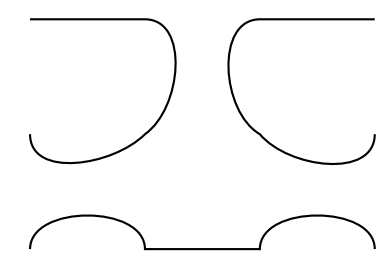

$+ \quad \sqrt{-t}(z_i - z_j)\sqrt{-t}(z_i - z_k)\sqrt{-t}(z_j - z_k)$ 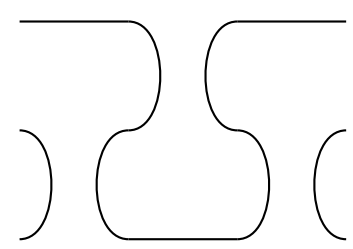 $+ \quad \sqrt{-t}(z_i - z_j)(tz_i + z_k)\sqrt{-t}(z_j - z_k)$ 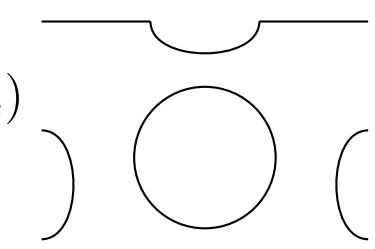

$+ \quad (tz_i + z_j)(tz_i + z_k)\sqrt{-t}(z_j - z_k)$ 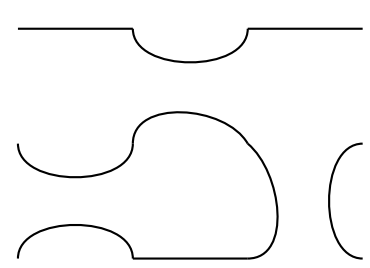 $+ \quad \sqrt{-t}(z_i - z_j)(tz_i + z_k)(tz_j + z_k)$ 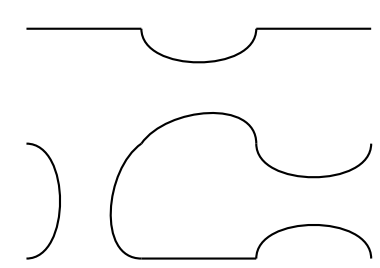

By Lemma 2.5 we can straighten curvy lines and substitute the loop in the sixth term with $(\sqrt{-t} + 1/\sqrt{-t})$. Thus, the eight terms in our expansion can be collected into five terms corresponding to the five crossingless matchings of six points on a circle. Indeed, notice the last four terms will gather into one single term:

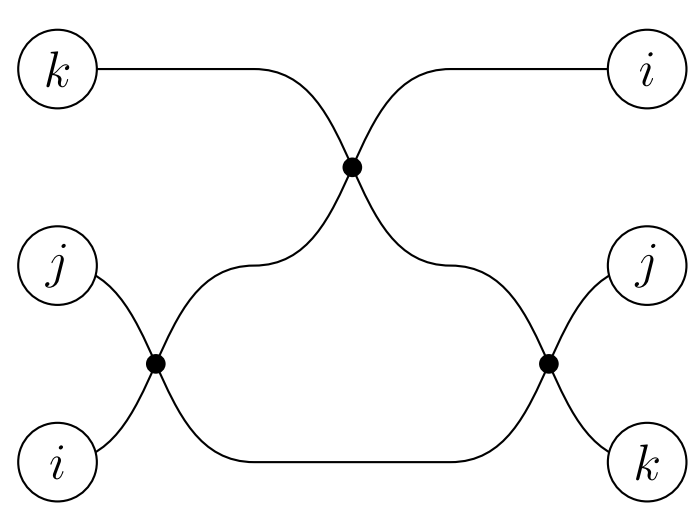

$=$

$= \quad (tz_i + z_j)(tz_i + z_k)(tz_j + z_k)$ 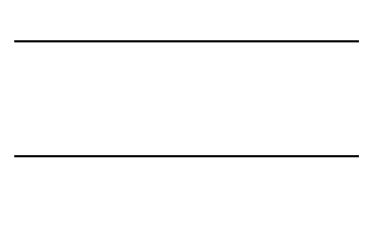 $+ \quad (tz_i + z_j)\sqrt{-t}(z_i - z_k)\sqrt{-t}(z_j - z_k)$ 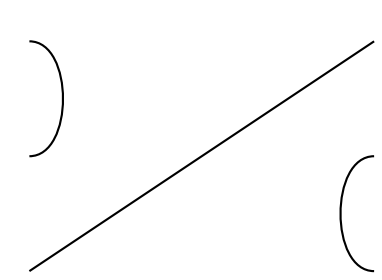

$$+ \quad \sqrt{-t}(z_i - z_j)\sqrt{-t}(z_i - z_k)(tz_j + z_k)$$ 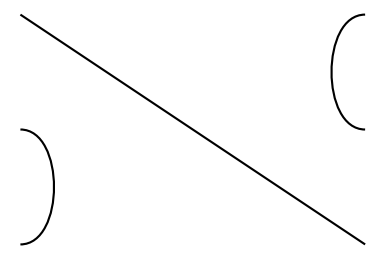 $$+ \quad (tz_i + z_j)\sqrt{-t}(z_i - z_k)(tz_j + z_k)$$ 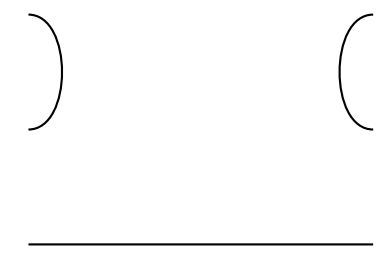

$$+ \left[\sqrt{-t}(z_i - z_j)\sqrt{-t}(z_i - z_k)\sqrt{-t}(z_j - z_k) + (\sqrt{-t} + \frac{1}{\sqrt{-t}})\sqrt{-t}(z_i - z_j)(tz_i + z_k)\sqrt{-t}(z_j - z_k) \right.$$
$$\left. + (tz_i + z_j)(tz_i + z_k)\sqrt{-t}(z_j - z_k) + \sqrt{-t}(z_i - z_j)(tz_i + z_k)(tz_j + z_k)\right]$$ 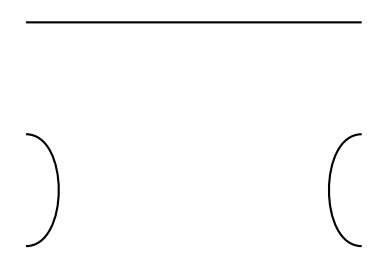

The complicated-looking coefficient we obtained in front of the fifth term simplifies quite a bit but we delay the simplification until after we have done the expansion for the right hand-side of the Yang-Baxter equation as well.

By skein relation (1) and Lemma 2.5 we have

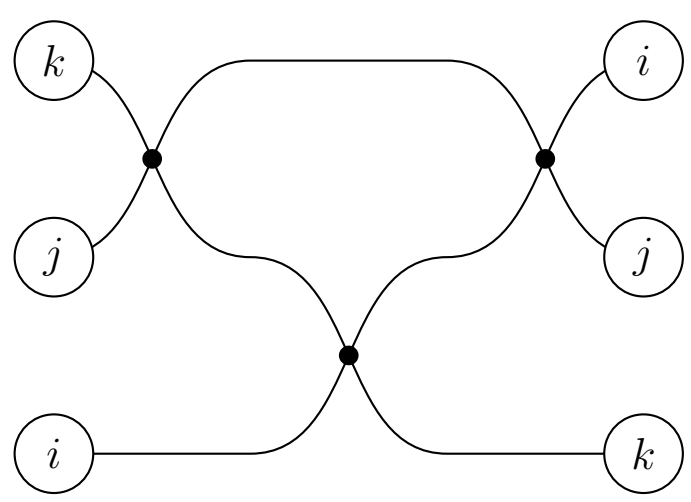

$=$

$$= \quad (tz_i + z_j)(tz_i + z_k)(tz_j + z_k)$$ 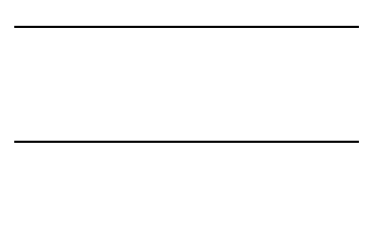 $$+ \quad (tz_i + z_j)\sqrt{-t}(z_i - z_k)\sqrt{-t}(z_j - z_k)$$ 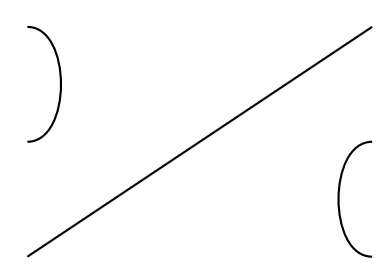

$$+ \quad \sqrt{-t}(z_i - z_j)\sqrt{-t}(z_i - z_k)(tz_j + z_k)$$ 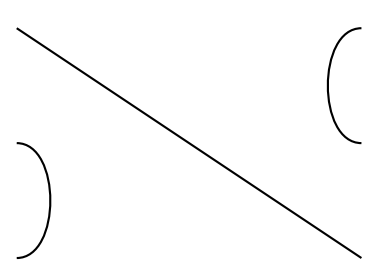 $$+ \quad (tz_i + z_j)\sqrt{-t}(z_i - z_k)(tz_j + z_k)$$ 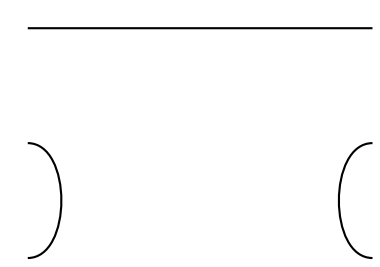

$$+ \left[\sqrt{-t}(z_i - z_j)\sqrt{-t}(z_i - z_k)\sqrt{-t}(z_j - z_k) + (\sqrt{-t} + \frac{1}{\sqrt{-t}})\sqrt{-t}(z_i - z_j)(tz_i + z_k)\sqrt{-t}(z_j - z_k) \right.$$
$$\left. + (tz_i + z_j)(tz_i + z_k)\sqrt{-t}(z_j - z_k) + \sqrt{-t}(z_i - z_j)(tz_i + z_k)(tz_j + z_k)\right]$$ 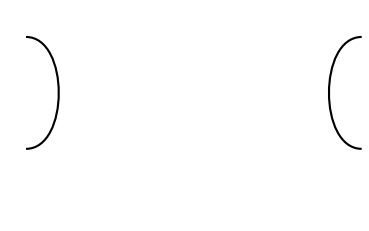

To prove the Yang-Baxter equation we have to prove that the two expansions we obtained are equal. The first three terms on both sides of the equation cancel with each other. The last two terms on the left-hand and right-hand side of the equation have their coefficients switched. Thus, the Yang-Baxter equation reduces to simply checking the equality

$$(tz_i + z_j)\sqrt{-t}(z_i - z_k)(tz_j + z_k) = \Big[\sqrt{-t}(z_i - z_j)\sqrt{-t}(z_i - z_k)\sqrt{-t}(z_j - z_k) + \\ + (\sqrt{-t} + \frac{1}{\sqrt{-t}})\sqrt{-t}(z_i - z_j)(tz_i + z_k)\sqrt{-t}(z_j - z_k) + \\ + (tz_i + z_j)(tz_i + z_k)\sqrt{-t}(z_j - z_k) + \\ + \sqrt{-t}(z_i - z_j)(tz_i + z_k)(tz_j + z_k)\Big] \tag{5}$$

which is not difficult to verify. Hence, the Boltzmann weights in Row 1 of Table 2 satisfy the Yang-Baxter equation.

**Second case:** We can repeat almost the same process for the Boltzmann weights in Row 2 of Table 2. For this reason we only provide an overview in this case.

Define the following cap weights:

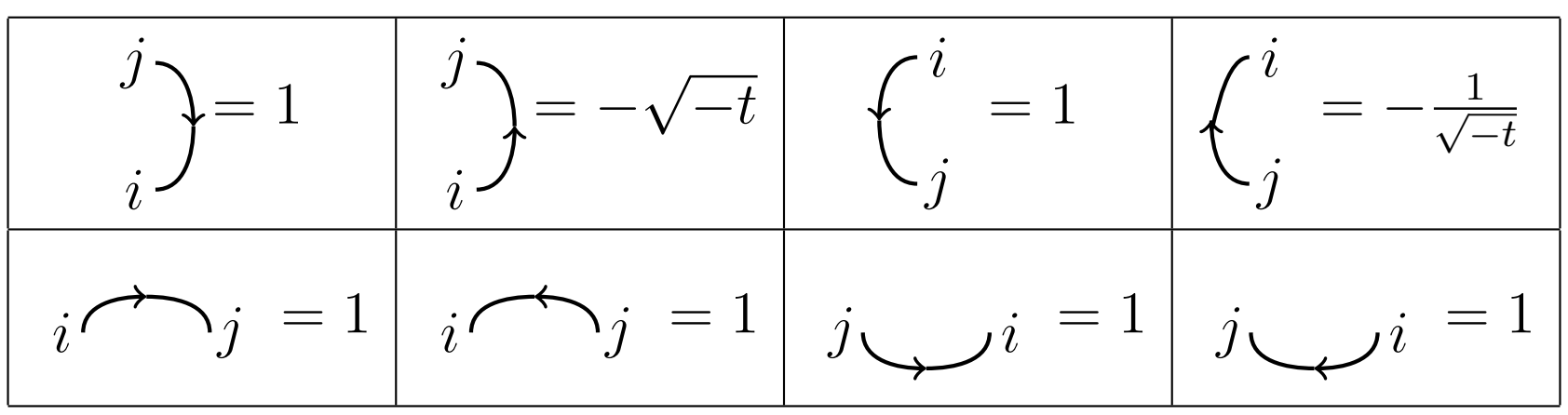

**Lemma 2.6.** *The Boltzmann weights in Row 2 of Table 2 satisfy the following skein relation:*

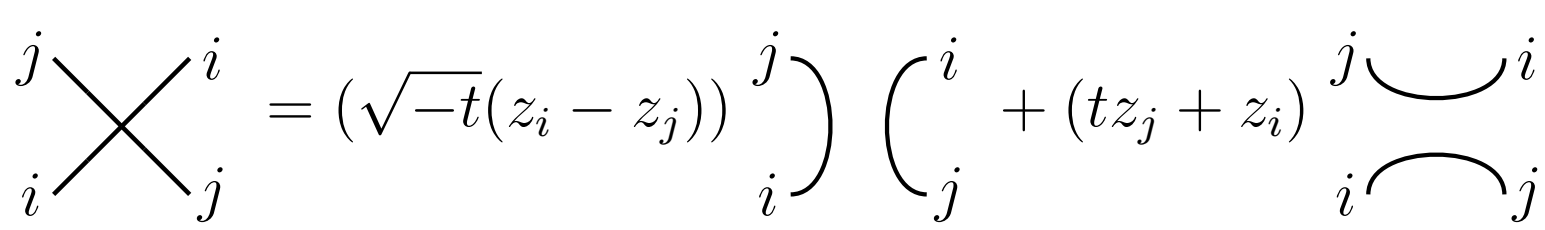

(6)

*Proof.* The proof follows by direct computation. □

The cap weights also satisfy all properties in Lemma 2.5 with the difference that the value of the loop changes to $-(\sqrt{-t} + 1/\sqrt{-t})$. Again one expands each side of the Yang-Baxter equation into five terms corresponding to the five crossingless matchings of six points on a circle.

As in the first case we end up with a single identity to be checked:

$$
\begin{aligned}
(tz_j + z_i)\sqrt{-t}(z_i - z_k)(tz_k + z_j) = \Big[ & \sqrt{-t}(z_i - z_j)\sqrt{-t}(z_i - z_k)\sqrt{-t}(z_j - z_k) - \\
& - \left(\sqrt{-t} + \frac{1}{\sqrt{-t}}\right)\sqrt{-t}(z_i - z_j)(tz_k + z_i)\sqrt{-t}(z_j - z_k) + \\
& + (tz_j + z_i)(tz_k + z_i)\sqrt{-t}(z_j - z_k) + \\
& + \sqrt{-t}(z_i - z_j)(tz_k + z_i)(tz_k + z_j) \Big]
\end{aligned} \tag{7}
$$

Remarkably, identity (7) is equivalent to identity (5) by just swapping the indices $i$ and $k$ followed by multiplying both sides of the equation by $-1$. This concludes the proof of Proposition 2.3. □

Lemma 2.2 and Proposition 2.3 are enough to give a proof of Theorem 2.1. The idea is to combine together the Boltzmann weights in Row 1 and Row 2 in Table 2 by using Lemma 2.2.

*Proof of Theorem 2.1.* By Lemma 2.2 all instances of the Yang-Baxter equation are divided into two groups. In the first we never have an $a_1$ state in the equation, in the second we never have an $a_2$ state in the equation. If we are in the first group we can "pretend" that our Boltzmann weights are the same as the Boltzmann weights in Row 2 in Table 2. Indeed, because we never have an $a_1$ state in these equation Table 1 and Row 2 in Table 2 become indistinguishable for the equations. Similarly, if we are in the second group, we use the Boltzmann weights in Row 1 in Table 2. Now Theorem 2.1 follows directly from Proposition 2.3. □

# 3 Application to Gamma-Gamma ice

In this section we show how we can apply the techniques developed in the proof of Theorem 2.1 to prove the Yang-Baxter equation for Gamma-Gamma ice. As mentioned above Gamma-Gamma ice is a free-fermionic six-vertex model introduced by Brubaker, Bump and Friedberg.

Gamma-Gamma ice uses the following Boltzmann weights:

| $a_1$ | $a_2$ | $b_1$ | $b_2$ | $c_1$ | $c_2$ |
|---|---|---|---|---|---|
| 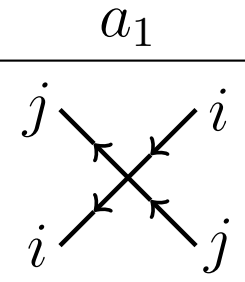 | 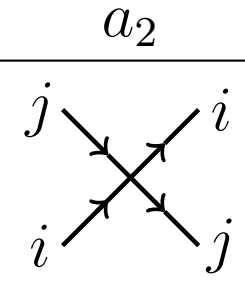 | 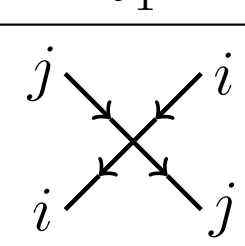 | 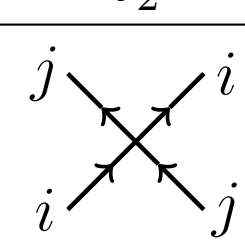 | 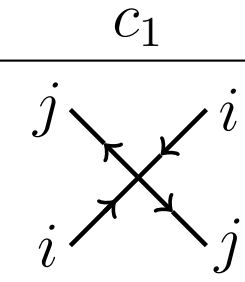 | 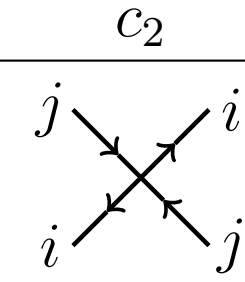 |
| $tz_i + z_j$ | $tz_j + z_i$ | $t(z_j - z_i)$ | $z_i - z_j$ | $(t+1)z_i$ | $(t+1)z_j$ |

Table 3: *Boltzmann weights for Gamma-Gamma ice.*

These Boltzmann weights satisfy the free-fermionic condition.

$$a_1a_2 + b_1b_2 = c_1c_2.$$

Note that the Boltzmann weights for Gamma-Gamma ice in Table 3 differ by the Boltzmann weights in Table 1 simply by a factor of $\sqrt{-t}$ that was moved between $b_1$ and $b_2$.

Brubaker et al. [BBF] prove that the weights in Table 3 satisfy the Yang-Baxter equation, however their argument relies on a computer calculation.

**Theorem 3.1** (Gamma-Gamma Yang-Baxter equation)**.** *For any way the external edges are fixed, the Boltzmann weights in Table 3 satisfy:*

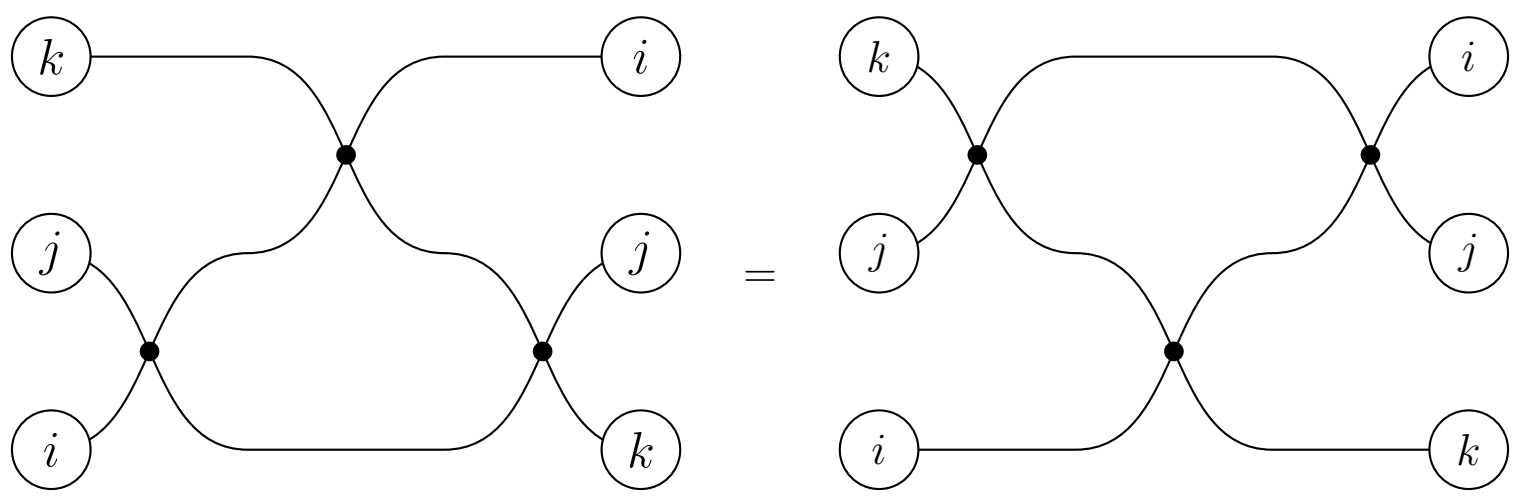

We provide a proof of Theorem 3.1 that utilises Proposition 2.3. But before giving the proof we need an important Lemma.

**Lemma 3.1.** *Let $a_1, a_2, b_1, b_2, c_1, c_2$ be the Boltzmann weights in Row 1 or Row 2 in Table 2. Then the Boltzmann weights $a_1, a_2, Cb_1, b_2/C, c_1, c_2$ satisfy the Yang-Baxter equation for any constant $C = f(t)$.*

*Proof.* We only prove the statement for the Boltzmann weights in Row 1 in Table 2 as the other case is analogous. We leverage the proof of Proposition 2.3. Modify the cap weights as follows:

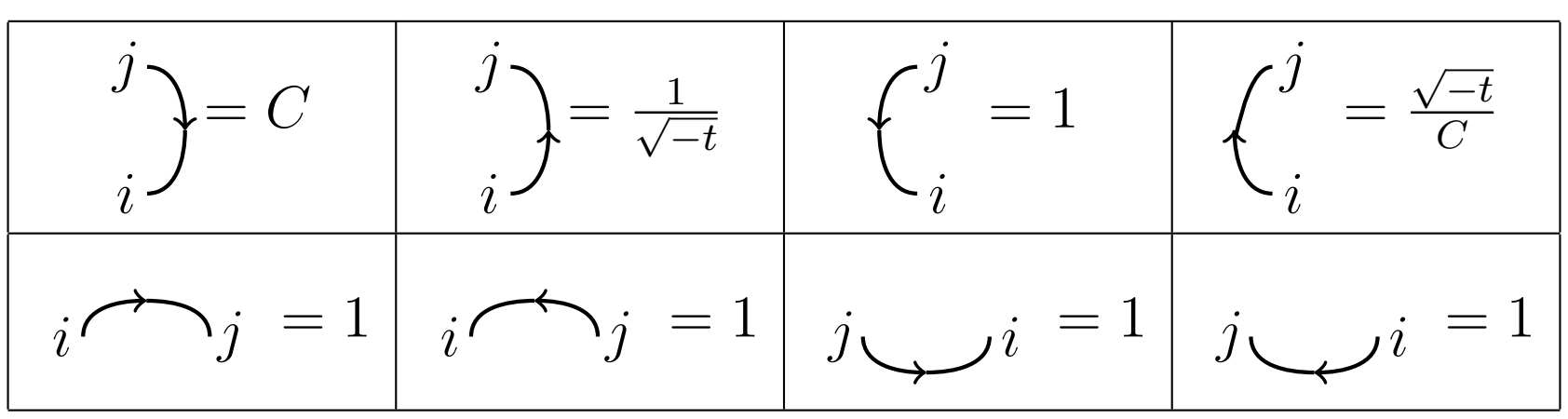

We use the same skein relation coefficients as in first case of the proof of Proposition 2.3.

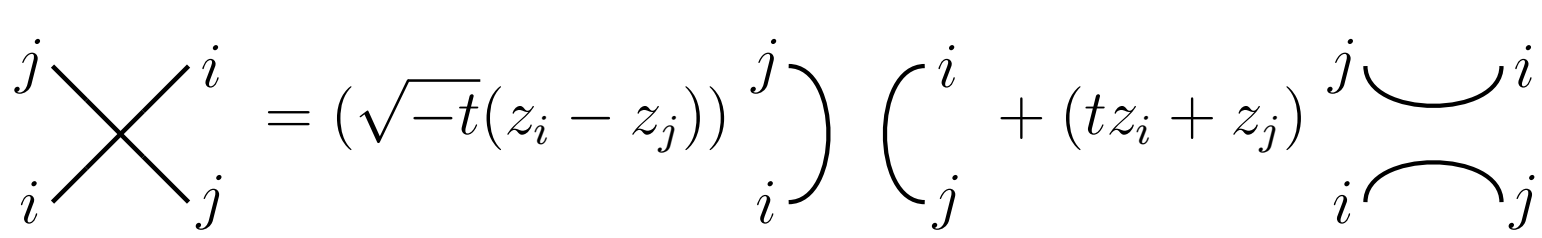

(8)

Notice that because we modified the cap weights, the skein relation gives the modified weights $a_1, a_2, Cb_1, b_2/C, c_1, c_2$. Nevertheless, the coefficients of the skein relation stayed the same.

Moreover, the cap weights satisfy almost the same properties as in Lemma 2.5.

**Lemma 3.2.** *The following identities hold:*

***Property 1:***

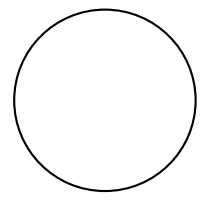

$$= \sqrt{-t} + \frac{1}{\sqrt{-t}} \tag{9}$$

*where we are summing over the two possible orientations of the loop.*

***Property 2:***

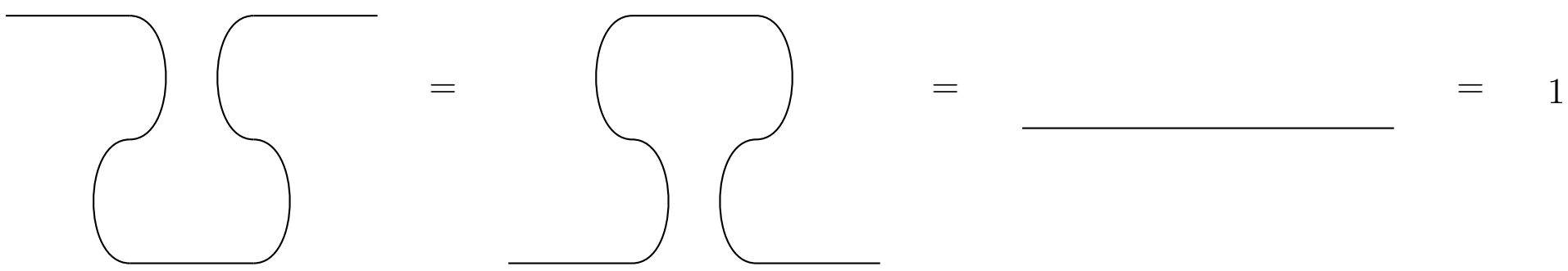

(10)

*where the equation holds for both orientations of the curvy lines and the straight line.*

One may expect that not having the second property of Lemma 2.5 could cause problems with the proof. But it turns out that this is not a problem. As before we could expand both sides of the Yang-Baxter equation into eight terms that can be collected into five terms corresponding to the five crossingless matchings of six points on a circle. The coefficients in front of all terms remain the same because the skein relation coefficients (8) and the loop in Lemma 3.2 all have the same value as in the proof of Proposition 2.3, i.e. the same as in skein relation (1) and the loop in Lemma 2.5. The only difference is that we cannot straighten the curvy line

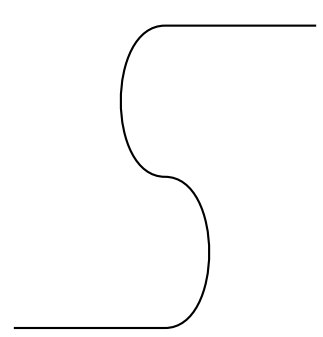

Thus, instead of having the terms

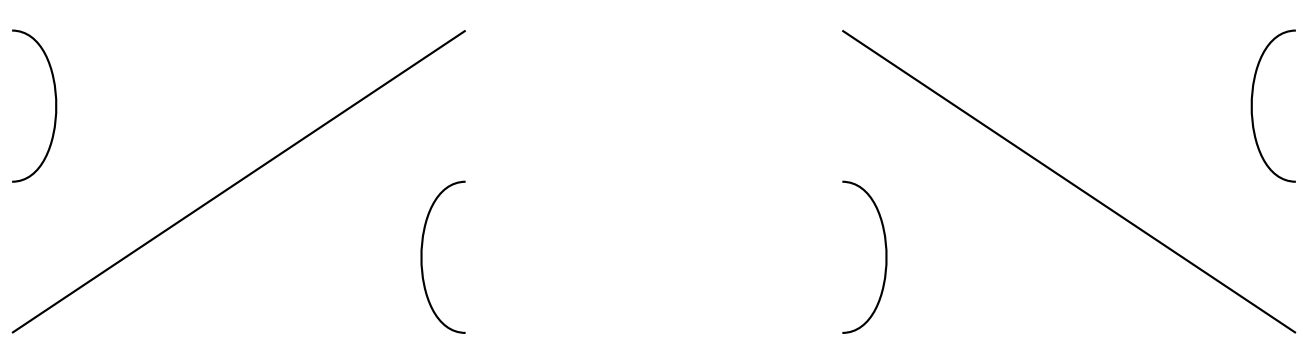

we have the terms

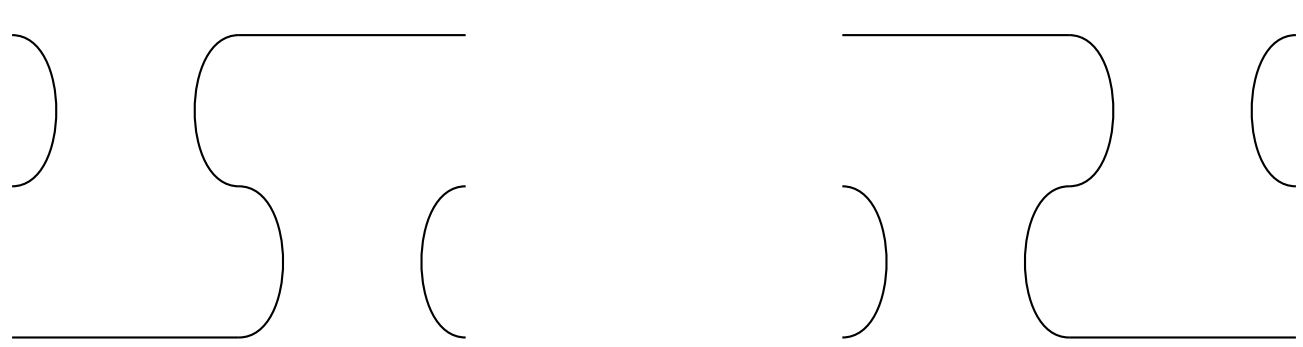

As a result, the curvy line itself also has a weight instead of just being equal to one. But this is not a problem as both sides of the expansion will have the same curvy line contributing the same weight in the respective terms. This is because the orientation of the curvy line is determined by the orientations of the external edges which are the same on both sides of the Yang-Baxter equation. Hence, the proof of Proposition 2.3 remains valid in this case as well. Thus, the Yang-Baxter equation once again follows from identity (5). This concludes the proof of Lemma 3.1. □

Lemma 3.1 allows us to prove Theorem 3.1.

*Proof of Theorem 3.1.* By setting $C = \sqrt{-t}$ in Lemma 3.1 we obtain that the following Boltzmann weights satisfy the Yang-Baxter equation:

| | $a_1$ | $a_2$ | $b_1$ | $b_2$ | $c_1$ | $c_2$ |
|---|---|---|---|---|---|---|
| | 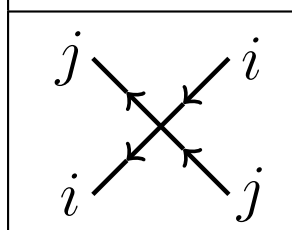 | 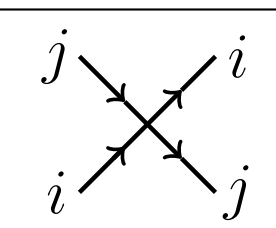 | 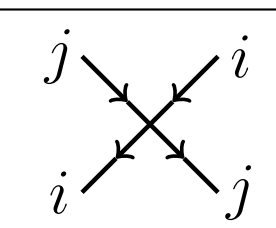 | 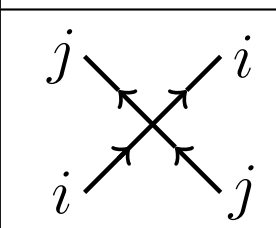 | 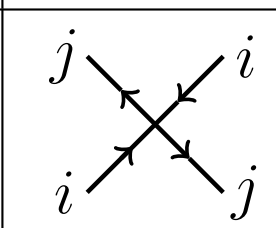 | 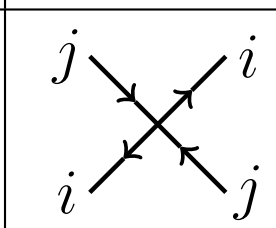 |
| Row 1 | $tz_i + z_j$ | $tz_i + z_j$ | $t(z_j - z_i)$ | $(z_i - z_j)$ | $(t+1)z_i$ | $(t+1)z_j$ |
| Row 2 | $tz_j + z_i$ | $tz_j + z_i$ | $t(z_j - z_i)$ | $(z_i - z_j)$ | $(t+1)z_i$ | $(t+1)z_j$ |

Table 4: *Modifying the weights in Table 2 by moving a constant between $b_1$ and $b_2$.*

Now by Lemma 2.2 we can mix the $a_1$ and $a_2$ weights and set $a_1 = tz_i + z_j$ and $a_2 = tz_j + z_i$. This concludes the proof. □

# 4 Kuperberg's field-free Boltzmann weights

This section is optional and can be skipped by the reader. It explains the connection between our methods and Kuperberg's [K]. The method we used to prove Proposition 2.3 comes from Kuperberg's proof of a Yang-Baxter equation for the field-free six vertex model.

In [K] Kuperberg uses a six-vertex model with Boltzmann weights that satisfy the field-free condition. The field-free condition means that $a_1 = a_2$ and $b_1 = b_2$. The condition $c_1 = c_2$ is optional.

Let $h$ be an arbitrary complex numbers. Now for any $x \in \mathbb{C}$ define $q = e^{hx}$. Furthermore, for convenience define $[x] = \frac{q^{x/2} - q^{-x/2}}{q^{1/2} - q^{-1/2}}$. Kuperberg's field-free weights are

| $a_1$ | $a_2$ | $b_1$ | $b_2$ | $c_1$ | $c_2$ |
|---|---|---|---|---|---|
| 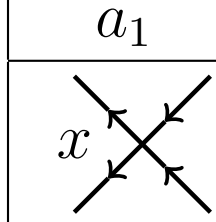 | 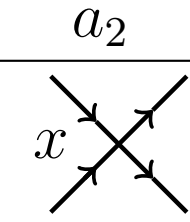 | 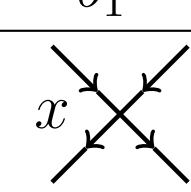 | 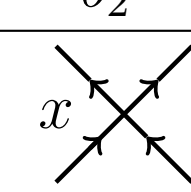 | 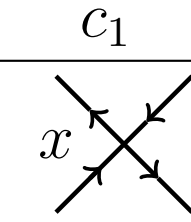 | 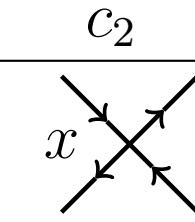 |
| $[x-1]$ | $[x-1]$ | $[x]$ | $[x]$ | $-q^{x/2}$ | $-q^{-x/2}$ |

Table 5: *Kuperberg's field-free weights.*

Kuperberg's weights are parameterised by the complex number $x$ and the value of $q$.

In [K] Kuperberg proves a Yang-Baxter equation by using skein-relations similarly to how we used them in Proposition 2.3.

**Theorem 4.1** (Kuperberg's field-free Yang-Baxter equation)**.** *Let $z = x + y$. For any way the external edges are fixed, the Boltzmann weights in Table 5 satisfy:*

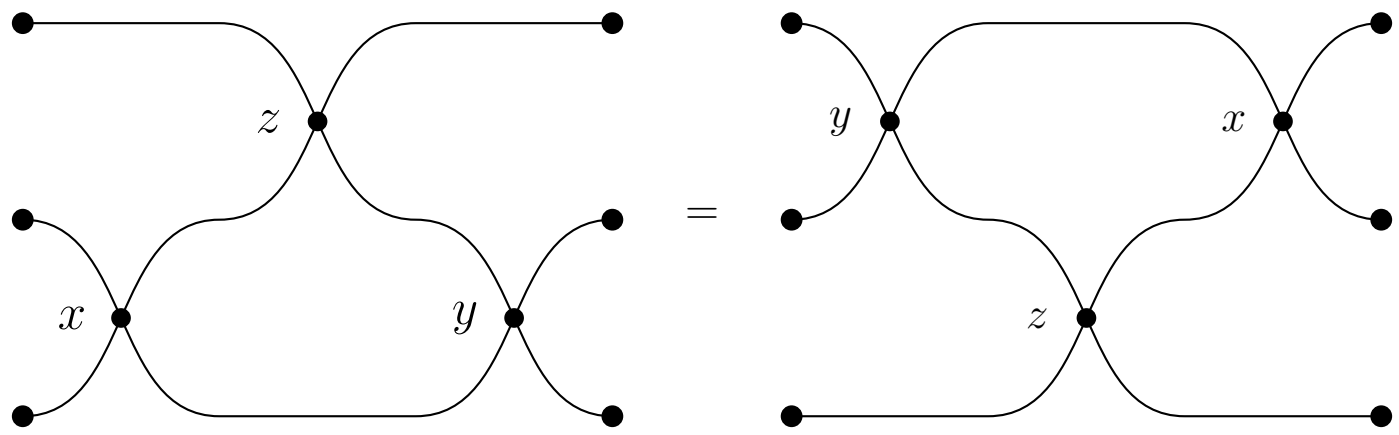

We do not present Kuperberg's proof of this fact as it is very similar to the proof of Proposition 2.3

Nevertheless, we do explain how our weights in Tables 1 and 2 are related to Kuperberg's weights in Table 5. We show that Kuperberg's weights in Table 5 can be substituted with the weights:

| $a_1$ | $a_2$ | $b_1$ |
|---|---|---|
| 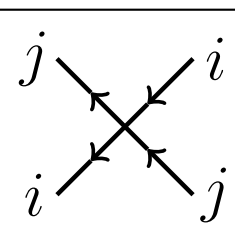 | 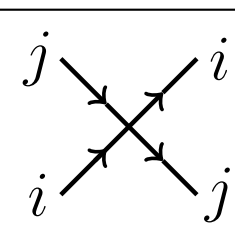 | 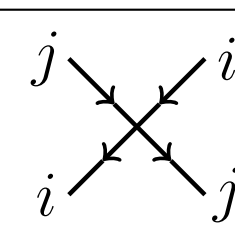 |
| $\dfrac{\frac{1}{\sqrt{-t}}\sqrt{\frac{z_i}{z_j}} - \sqrt{-t}\sqrt{\frac{z_j}{z_i}}}{\sqrt{-t} - \frac{1}{\sqrt{-t}}}$ | $\dfrac{\frac{1}{\sqrt{-t}}\sqrt{\frac{z_i}{z_j}} - \sqrt{-t}\sqrt{\frac{z_j}{z_i}}}{\sqrt{-t} - \frac{1}{\sqrt{-t}}}$ | $\dfrac{\sqrt{\frac{z_i}{z_j}} - \sqrt{\frac{z_j}{z_i}}}{\sqrt{-t} - \frac{1}{\sqrt{-t}}}$ |
| $b_2$ | $c_1$ | $c_2$ |
| 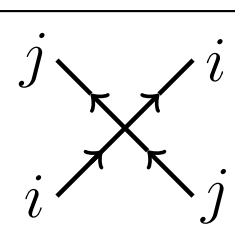 | 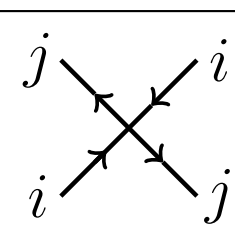 | 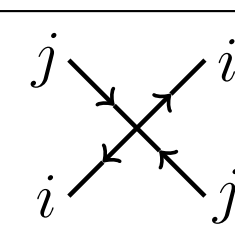 |
| $\dfrac{\sqrt{\frac{z_i}{z_j}} - \sqrt{\frac{z_j}{z_i}}}{\sqrt{-t} - \frac{1}{\sqrt{-t}}}$ | $-\sqrt{\dfrac{z_i}{z_j}}$ | $-\sqrt{\dfrac{z_j}{z_i}}$ |

Table 6: *The Table 5 Boltzmann weights after the substitutions.*

It is worth noting that we just substituted the $x, y, z$ labels of the vertices with the row labels $i, j, k$. The reasons this happens is that when transforming the weights, the variable $x$ will be a function of $i, j$, the variable $y$ will be a function of $j, k$ and the variable $z$ will be a function of $i, k$.

Indeed, let $q^{1/2} = \sqrt{-t}$, $q^x = \frac{z_i}{z_j}$, $q^y = \frac{z_j}{z_k}$ and $q^z = \frac{z_i}{z_k}$. It is not difficult to show this is a valid substitution. One can check that the weights in Table 5 become the weights in Table 6 after the substitutions are made.

As a corollary we get that the new weights also satisfy the Yang-Baxter equation:

**Theorem 4.2.** *For any way the external edges are fixed, the Boltzmann weights in Table 6 satisfy:*

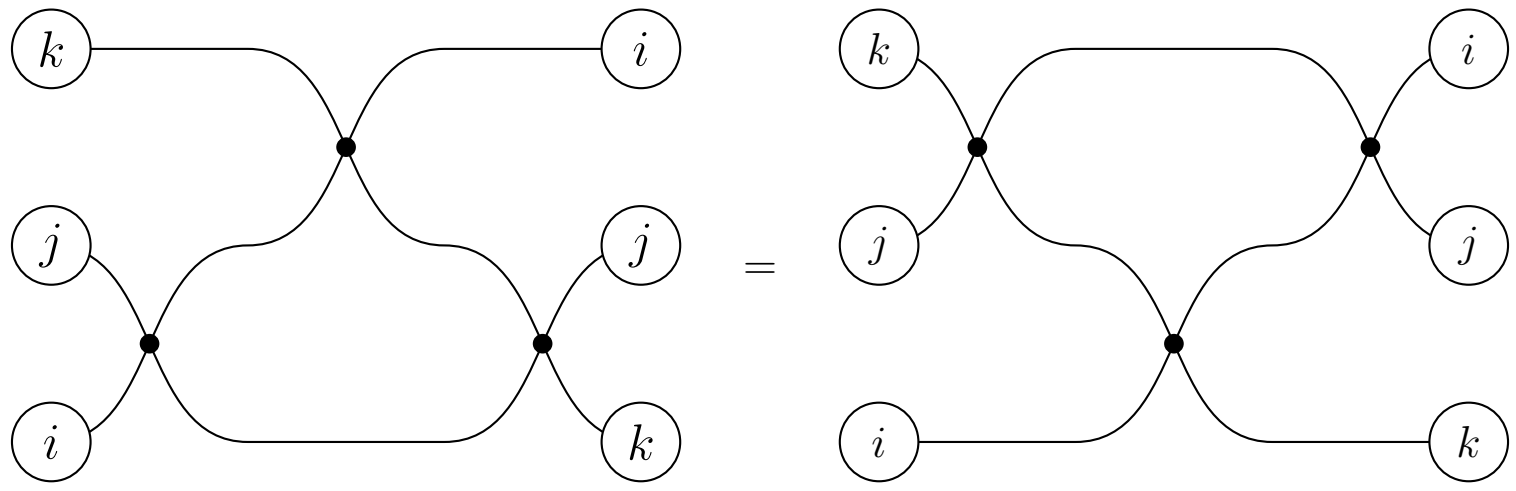

*Proof.* Observe that the condition $x + y = z$ is equivalent to the condition $q^x q^y = q^z$. This is true as

$$q^x q^y = \frac{z_i}{z_j}\frac{z_j}{z_k} = \frac{z_i}{z_k} = q^z$$

Thus, the theorem follows directly from Theorem 4.1. □

Theorem 4.2 allows us to relate Kuperberg's results to Table 1 and Table 2.

| | $a_1$ | $a_2$ | $b_1$ | $b_2$ | $c_1$ | $c_2$ |
|---|---|---|---|---|---|---|
| | 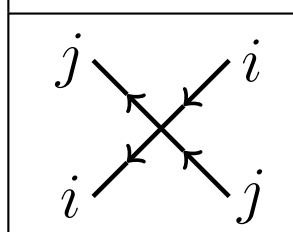 | 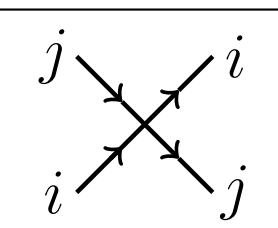 | 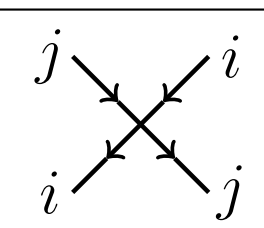 | 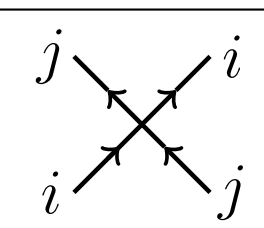 | 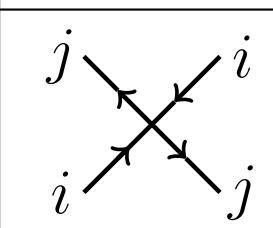 | 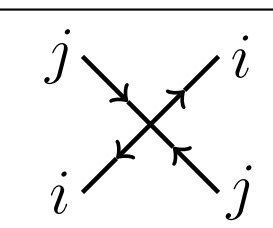 |
| Row 1 | $tz_j + z_i$ | $tz_j + z_i$ | $\sqrt{-t}(z_i - z_j)$ | $\sqrt{-t}(z_i - z_j)$ | $(t+1)z_i$ | $(t+1)z_j$ |
| Row 2 | $tz_i + z_j$ | $tz_i + z_j$ | $\sqrt{-t}(z_j - z_i)$ | $\sqrt{-t}(z_j - z_i)$ | $(t+1)z_i$ | $(t+1)z_j$ |

Table 7: *Two different transformations of Kuperberg's field-free weights.*

**Corollary 4.1.** *The Boltzmann weights in both Row 1 and Row 2 in Table 7 satisfy the Yang-Baxter equation.*

*Proof.* Note that each row of Boltzmann weights satisfying the Yang-Baxter equation can be multiplied by a symmetric function of $i$ and $j$ and the new row will still satisfy the Yang-Baxter equation.

To get the Boltzmann weights in Row 1 in Table 7 multiply the Boltzmann weights in Table 6 by $\sqrt{-t}(\sqrt{-t} - \frac{1}{\sqrt{-t}})\sqrt{z_i z_j}$.

To get the Boltzmann weights in Row 2 in Table 7 multiply the Boltzmann weights in Table 6 by $-\frac{1}{\sqrt{-t}}(\sqrt{-t} - \frac{1}{\sqrt{-t}})\sqrt{z_i z_j}$. Then abusing notation set $\frac{1}{\sqrt{-t}} = \sqrt{-t}$. □

The weights in Row 2 in Table 7 differ from the weights in Row 2 in Table 2 by a minus sign between $b_1$ and $b_2$. Nevertheless, due to Lemma 3.1, we see that Kuperberg's field-free weights are equivalent to the field-free weights in Table 2.

This gives an explanation as to why Kuperberg's methods in [K] work in our specific setting despite the fact that the weights in Table 1 are free-fermionic.

# 5 General form of Gamma-Gamma ice

In their article [BBF] Brubaker et al. prove a more general identity than Theorem 3.1 by using the weights in Table 8.

| $a_1$ | $a_2$ | $b_1$ | $b_2$ | $c_1$ | $c_2$ |
|---|---|---|---|---|---|
| 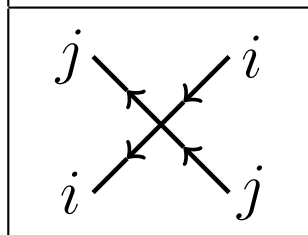 | 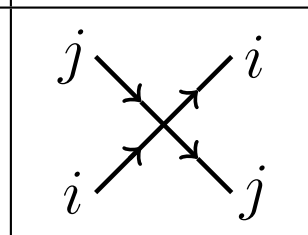 | 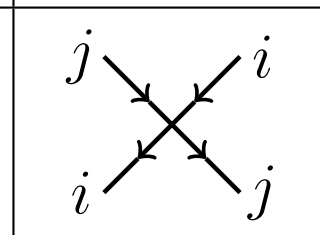 | 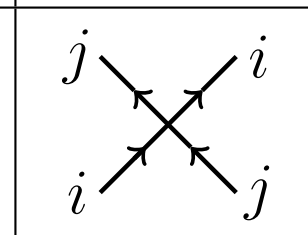 | 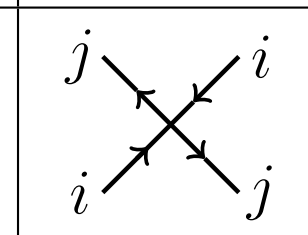 | 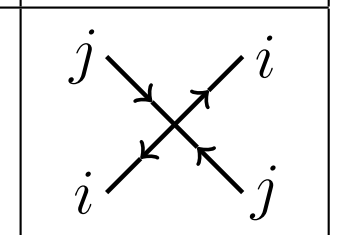 |
| $t_j z_i + z_j$ | $t_i z_j + z_i$ | $t_i z_j - t_j z_i$ | $z_i - z_j$ | $(t_i + 1) z_i$ | $(t_j + 1) z_j$ |

Table 8: *Generalised Boltzmann weights for Gamma-Gamma ice.*

Notice that the weights in Table 8 are exactly the weights in Table 3 when we take $t_i = t_j$.

**Theorem 5.1.** *The weights in Table 8 satisfy the Yang-Baxter equation*

An important questions is whether it is possible to use skein relations to prove Theorem 5.1. One might wonder if our approach in Sections 2 and 3 will work. The answer is yes but one has to significantly change the skein relation and the proof overall is more complicated. The reason for the complications is that $b_1$ and $b_2$ no longer have the factor $z_i - z_j$ in common. Moreover, the fact that $t_i \neq t_j$ causes other complications as well.

Thus, one may wonder how to go about discovering a skein relation in this case. What values of the cap weights would be suitable? We answer all these questions in the next two sections.

# 6 The Fish equation, cap weights and skein relations

In this section we study the fish equation in order to find suitable values for the cap weights. While the Yang-Baxter equation can be interpreted as a Reidemeister Type 3 move, the fish equation can be interpreted as a Reidemeister Type 1 move.

In our setting the fish equation is the identity

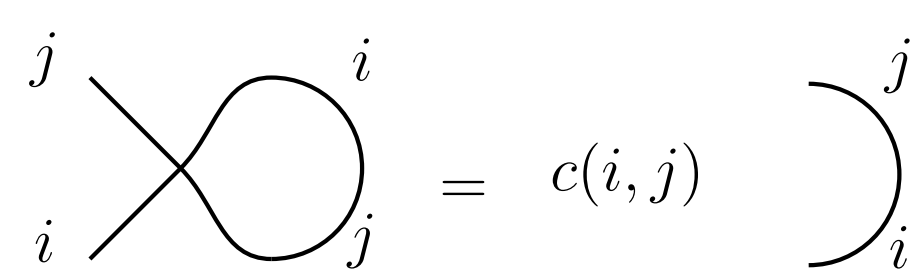

where $c(i, j)$ is a constant possibly dependent on $z_i, z_j, t_i, t_j$.

The fish equation is equivalent to the following two equations:

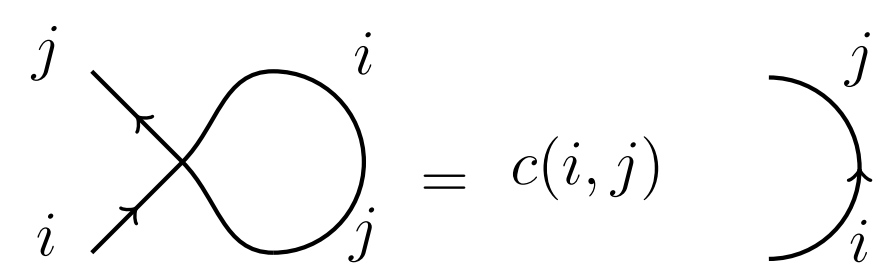

(11)

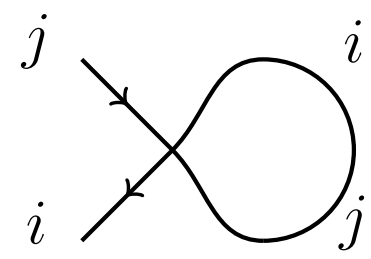 $=\ c(i,j)$ 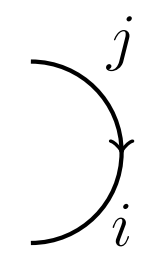 (12)

In both equations (11) and (12) we sum over the two possible states on the left-hand side. Each of the two possible states is just a vertex times a cap.

We translate the fish equation into an algebraic form. Let us denote the cap weights by

| | |
|---|---|
| 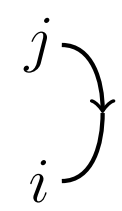 $= \alpha(i,j)$ | 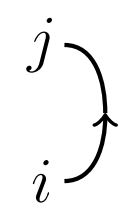 $= \beta(i,j)$ |

The algebraic forms of equations (11) and (12) are respectively

$$(t_i z_i + z_i) \cdot \beta(j,i) + (z_i - z_j) \cdot \alpha(j,i) = c(i,j) \cdot \beta(i,j) \tag{13}$$

$$(t_j z_j + z_j) \cdot \alpha(j,i) + (t_i z_j - t_j z_i) \cdot \beta(j,i) = c(i,j) \cdot \alpha(i,j) \tag{14}$$

Note that on the left-hand side the $i$ and $j$ are switched in the functions $\alpha$ and $\beta$. This represents the fact that the "head" of the fish interchanges $i$ and $j$.

**Theorem 6.1.** *We have the following solutions to the right fish equation:*

| $c(i,j)$ | $\alpha(i,j)$ | $\beta(i,j)$ |
|---|---|---|
| $a_2 = t_i z_j + z_i$ | $t_j$ | $-1$ |
| $a_2 = t_i z_j + z_i$ | $z_j$ | $z_i$ |
| $a_1 = t_j z_i + z_j$ | $(t_j+1)z_j t_i$ | $(t_i+1)z_i$ |
| $a_1 = t_j z_i + z_j$ | $-(t_j+1)$ | $(t_i+1)$ |

*Moreover, any linear combination (with scalars independent from $i,j$) of solutions which share the same constant is also a solution.*

*Proof.* The first statement is shown by direct computation. We provide more intuition about how we arrived at these solutions in the appendix. The second statement follows form the fact that $\alpha$ and $\beta$ have degree 1 on both sides of the equation in both (13) and (14). □

Hence, we have found possible values for the right caps.

**Corollary 6.1.** *Suppose $t_i = t_j = t$ in the right fish equation. Then we have the following solutions*

| $c(i,j)$ | $\alpha(i,j)$ | $\beta(i,j)$ |
|---|---|---|
| $a_2 = tz_j + z_i$ | $t$ | $-1$ |
| $a_2 = tz_j + z_i$ | $z_j$ | $z_i$ |
| $a_1 = tz_i + z_j$ | $z_j t$ | $z_i$ |
| $a_1 = tz_i + z_j$ | $-1$ | $1$ |

*Moreover, any linear combination (with scalars independent from $i, j$) of solutions which share the same constant is also a solution.*

*Proof.* Set $t_i = t_j = t$ in the solutions from Theorem 6.1. In the third and fourth row of solutions we remove the factor $(t+1)$ as it is now independent from $i$ and $j$. $\square$

**Remark 6.2.** *Note that the cap weights $\alpha = \sqrt{-t}$ and $\beta = 1/\sqrt{-t}$ from the proof of Theorem 3.1 in Section 3 differ by a factor of $-\sqrt{-t}$ from the first row of solutions in Corollary 6.1. Thus, one way we could have found the skein relation in section 3 would have been to first solve the fish equation and then get suitable values for the cap weights.*

We repeat the same process with the left fish equation in order to obtain values for the left caps.

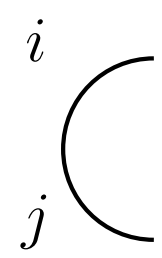

$c(i,j) \quad =$

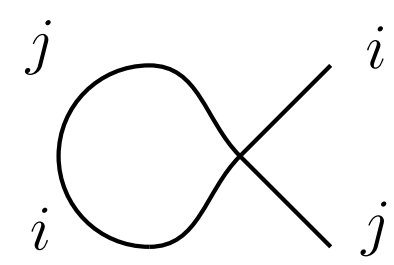

where abusing notation we again used $c(i,j)$ to denote the constant.

Denote the caps by

| | |
|---|---|
| 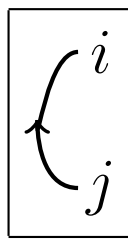 $= \gamma(i,j)$ | 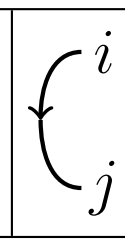 $= \delta(i,j)$ |

The algebraic forms of the two fish equations are

$$\gamma(i,j) \cdot c(i,j) = \gamma(j,i) \cdot (t_j z_j + z_j) + \delta(j,i) \cdot (z_i - z_j) \tag{15}$$

$$\delta(i,j) \cdot c(i,j) = \delta(j,i) \cdot (t_i z_i + z_i) + \gamma(j,i) \cdot (t_i z_j - t_j z_i) \tag{16}$$

We again find solutions to these equations.

**Theorem 6.2.** *We have the following solutions to the left fish equation:*

| $c(i,j)$ | $\gamma(i,j)$ | $\delta(i,j)$ |
|---|---|---|
| $a_1 = t_j z_i + z_j$ | $1$ | $t_i$ |
| $a_1 = t_j z_i + z_j$ | $z_j$ | $-z_i$ |
| $a_2 = t_i z_j + z_i$ | $-(t_j+1)z_j$ | $(t_i+1)z_i t_j$ |
| $a_2 = t_i z_j + z_i$ | $(t_j+1)$ | $(t_i+1)$ |

*Moreover, any linear combination (with scalars independent from $i, j$) of solutions which share the same constant is also a solution.*

*Proof.* The theorem follows by direct computation. We provide some intuition about how we found these solutions in the appendix $\square$

Hence, we have found possible values for the left cap weights as well.

**Corollary 6.3.** *Suppose $t_i = t_j = t$ in the right fish equation. Then we have the following solutions*

| $c(i,j)$ | $\gamma(i,j)$ | $\delta(i,j)$ |
|---|---|---|
| $a_1 = tz_i + z_j$ | $1$ | $t$ |
| $a_1 = tz_i + z_j$ | $z_j$ | $-z_i$ |
| $a_2 = tz_j + z_i$ | $-z_j$ | $z_i t$ |
| $a_2 = tz_j + z_i$ | $1$ | $1$ |

*Moreover, any linear combination (with scalars independent from $i, j$) of solutions which share the same constant is also a solution.*

*Proof.* Set $t_i = t_j = t$ in the solutions from Theorem 6.1. In the third and fourth row of solutions we remove the factor $(t + 1)$ as it is now independent from $i$ and $j$. □

**Remark 6.4.** *Note that the cap weights $\gamma = 1$ and $\delta = 1$ from the proof of Theorem 3.1 in Section 3 are the same as the solutions in the fourth row of Corollary 6.3. Thus, one way we could have found the skein relation in Section 3 would have been to first solve the fish equation and get suitable values for the cap weights.*

We present the following questions as interesting problems for further research:

**Question 6.5.** *Categorise all possible solutions of the right and left fish equations.*

**Question 6.6.** *What are the possible applications to solvable lattice models or knot theory for each of the solutions to the fish equation?*

# 7 Proof of general form of Gamma-Gamma ice

In this section we use skein relations to prove the general form of the Gamma-Gamma ice Yang-Baxter equation in Theorem 5.1. The cap weights that we use in this section are taken from our solutions to the fish equation.

Lemma 2.2 implies that we can divide the instances of the Yang-Baxter equation into two groups. In the first group we have at least four external edges pointing to the left and we know there is no $a_2$ state present. In the second group there are at least four external edges pointing to the right and there is no $a_1$ state. This allows us to use an approach similar to the proof of Theorem 2.1 where we use two different skein relations depending on which of the two states $a_1$ and $a_2$ is present. Furthermore, we have additional information about the orientation of the external edges in each case.

*Proof of Theorem 5.1.* We divide the proof in two cases.

**First case:** Let there be at least four external edges pointing to the right. By Lemma 2.2 we have that there is no $a_1$ state in any instance of the Yang-Baxter equation. Hence, we may pretend that we are working with only five weights in this case: $a_2$, $b_1$, $b_2$, $c_1$ and $c_2$. Our skein relation will be defined only for these five vertices.

Define the following cap weights and diagonal line weights:

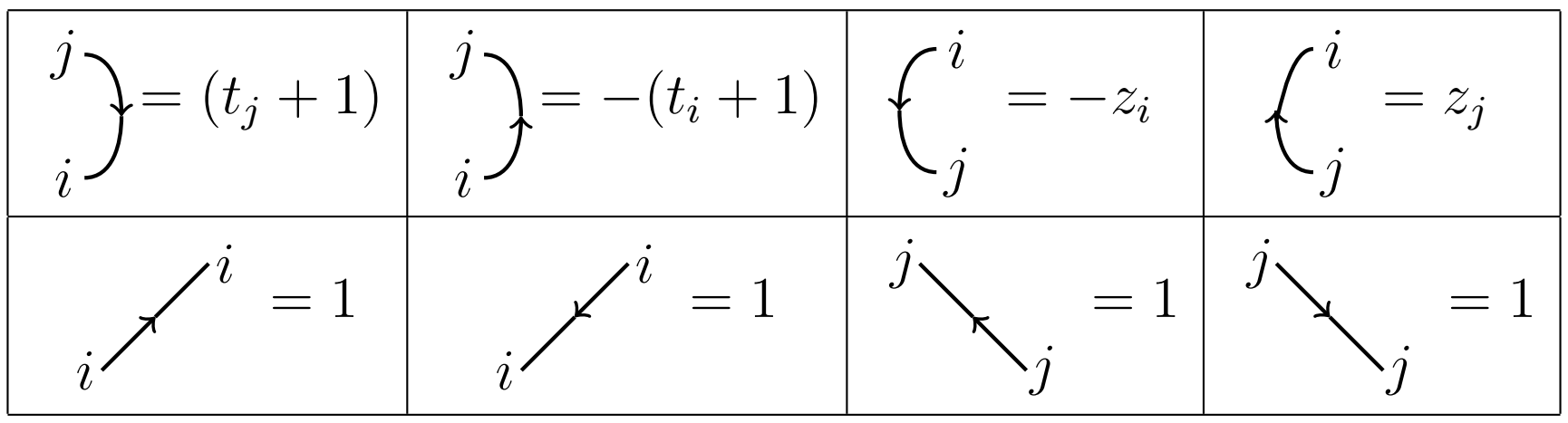

Notice that the cap weights are coming from our solutions to the fish equation in Theorems 6.1 and 6.2.

These cap and line weights lead to three important lemmas.

**Lemma 7.1.** *The Boltzmann weights in Table 8, excluding $a_1$, satisfy the following skein relation*

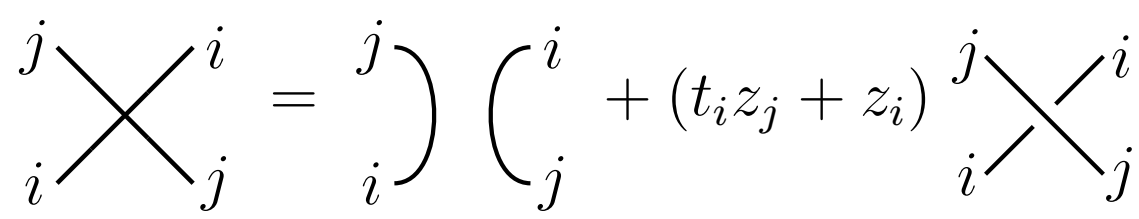

(17)

*where the second term in the skein relation represents the crossing of two straight diagonal lines where the upper line crosses above the lower line.*

*Proof.* The proof follows by direct computation. □

The fact that the diagonal lines contribute a weight of 1 gives us the following property that resembles the Reidemeister Type 2 move and is straightforward to check.

**Lemma 7.2.** *The following identity holds*

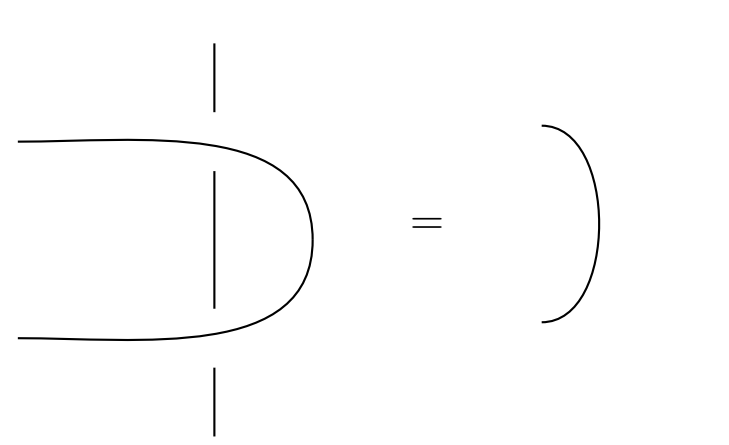

*where the straight line on the right-hand side contributes a factor of* 1.

Finally, we have the following intriguing identities which come up later in our proof.

**Lemma 7.3.** *The following identities hold*

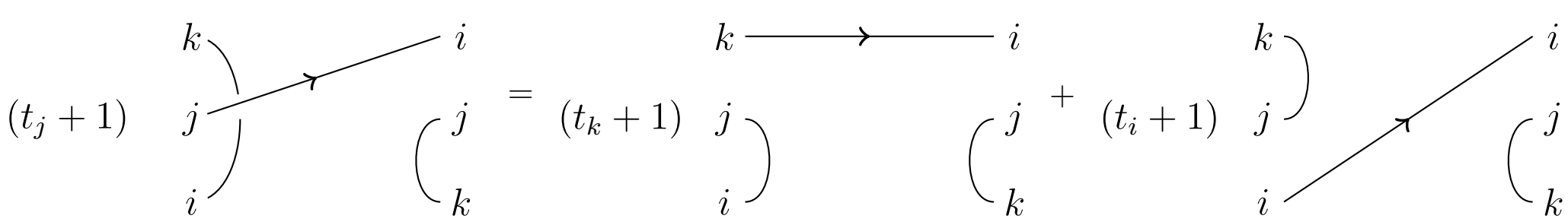

(18)

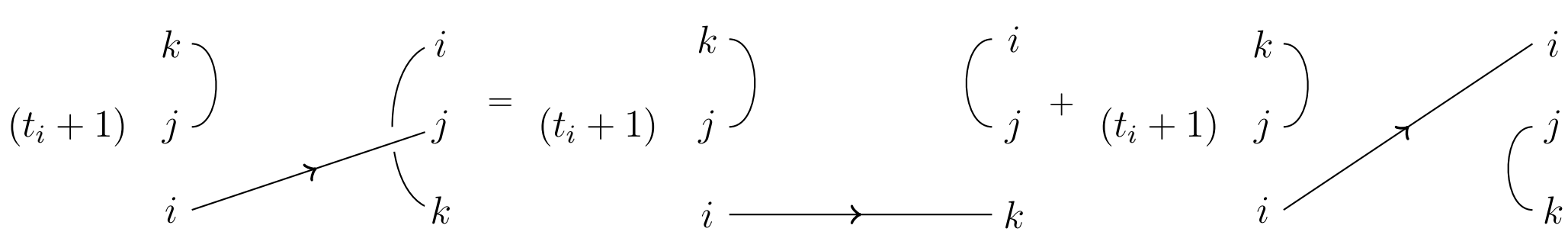

(19)

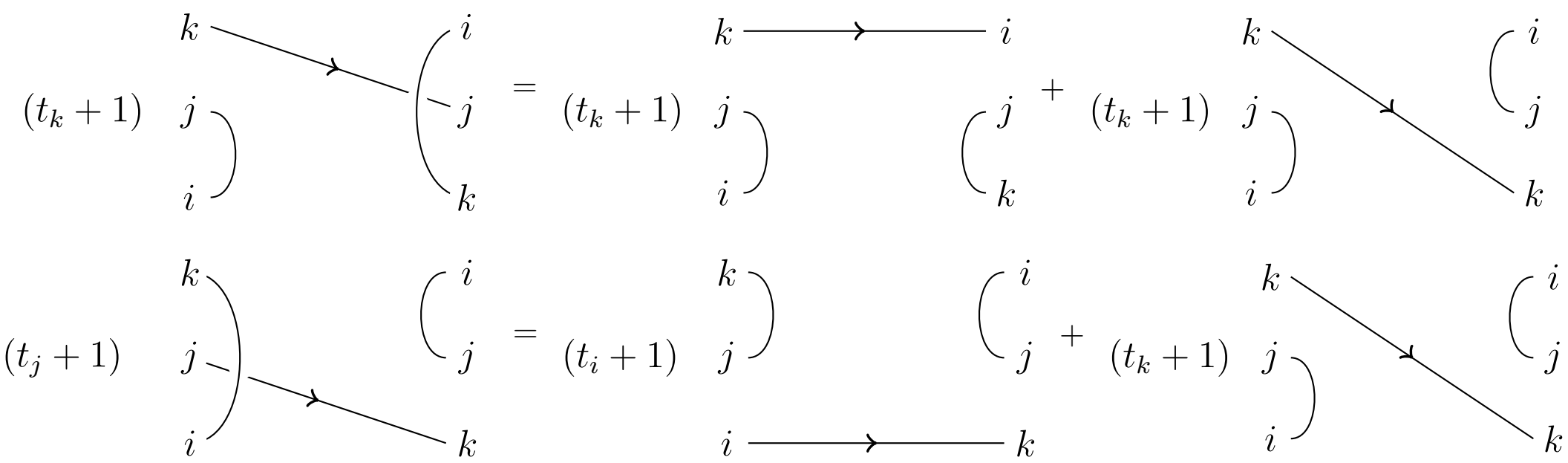

(20)

(21)

*where the big caps connecting the rows $i$ and $k$ should be considered as ordinary caps.*

*Proof.* Notice that identities (18) and (21) are equivalent. The same is true for identities (19) and (20). Now identities (18) and (19) are straightforward to check by direct computation. □

**Remark 7.4.** *The four identities in Lemma 7.3 resemble skein relations. Indeed, notice that if we substitute the crossing on the left-hand side with two vertical caps or two horizontal caps, respectively, we get the shapes on the right-hand side.*

A direct consequence of Lemma 7.3 is the next corollary:

**Corollary 7.5.** *The following triple equality holds.*

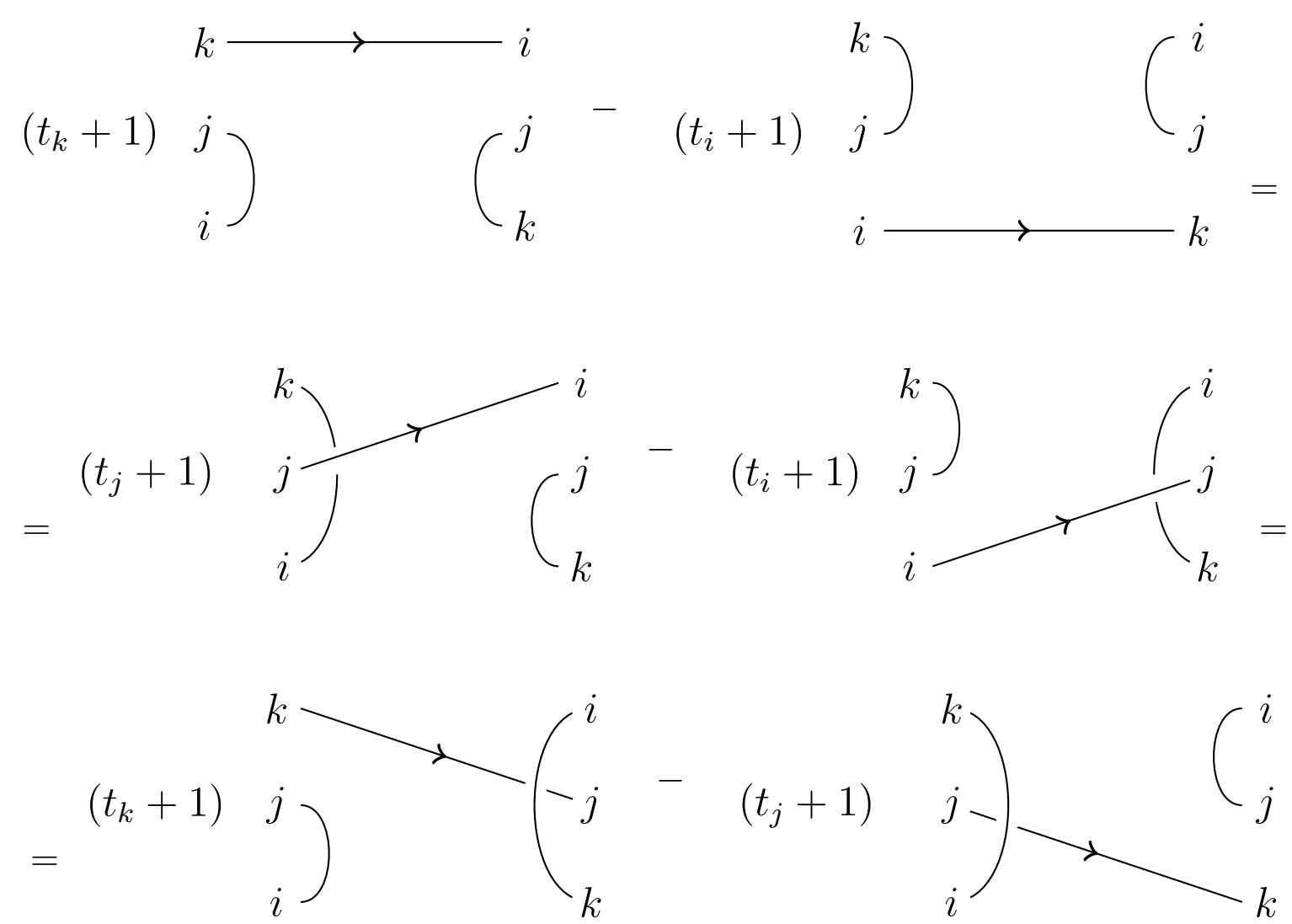

(22)

Now we are ready to prove the Yang-Baxter equation when there are at least four external edges pointing to the right. By using the skein relation (17) we can expand the left-hand side of the Yang-Baxter equation into eight terms.

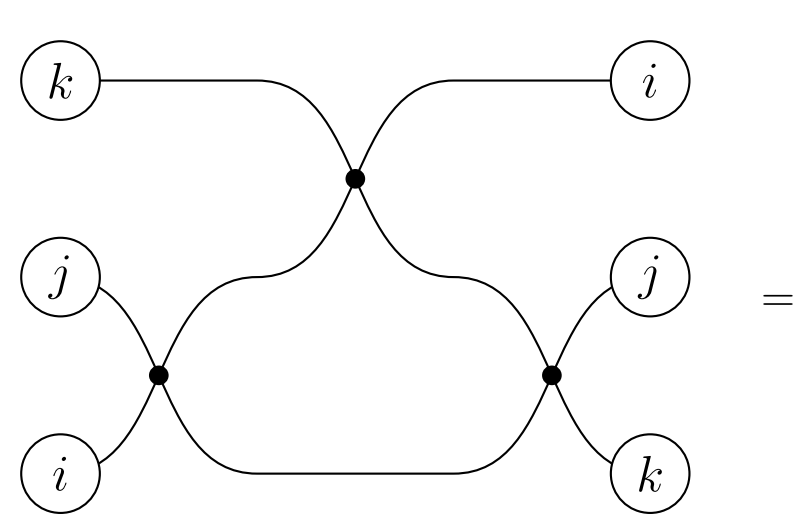

$=$ $(t_i z_j + z_i)(t_i z_k + z_i)(t_j z_k + z_j)$ 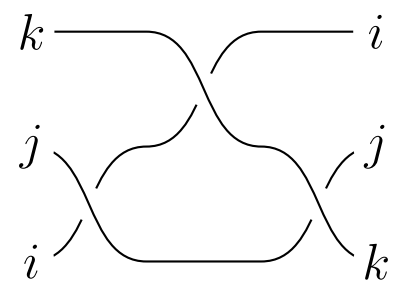

$+$ $(t_i z_j + z_i)(t_i z_k + z_i)$ 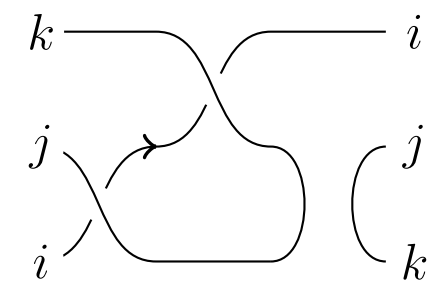

$+$ $(t_i z_k + z_i)(t_j z_k + z_j)$ 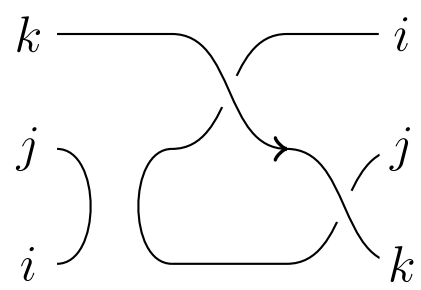

$+$ $(t_i z_j + z_i)(t_j z_k + z_j)$ 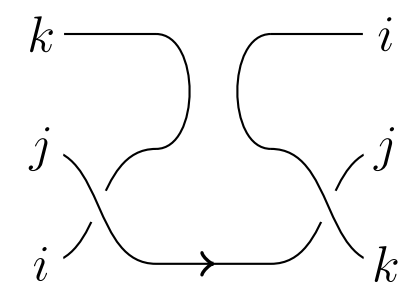

$+$ $(t_i z_k + z_i)$ 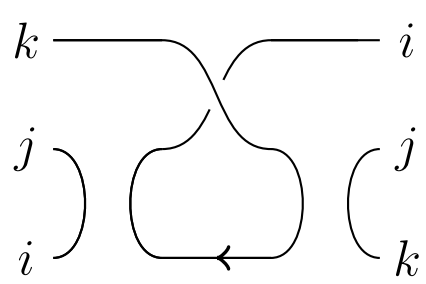

$+$ $(t_i z_j + z_i)$ 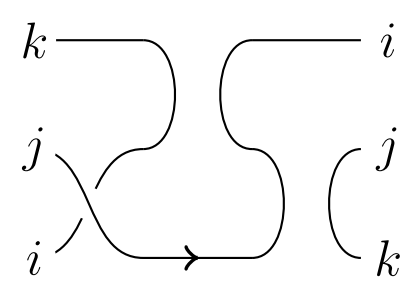

$+$ $(t_j z_k + z_j)$ 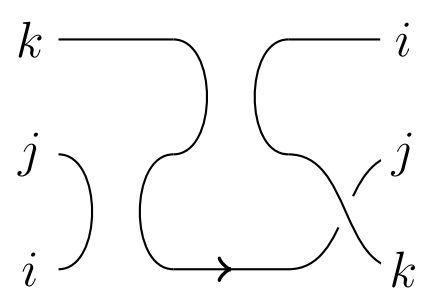

$+$ 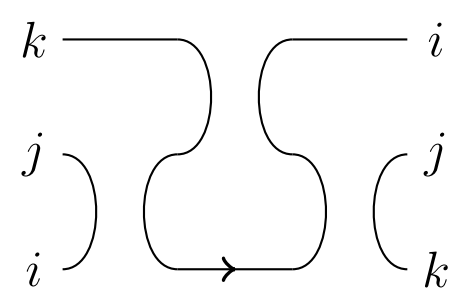

Notice that in all terms, except the first one, we have fixed the direction of one of the lines. Because we are in the case with at least four external edges pointing to the right, only terms with this orientation contribute any weight.

One can apply Lemma 7.2 to the second and third term but before we do this we repeat the same process for the right-hand side of the Yang-Baxter equation.

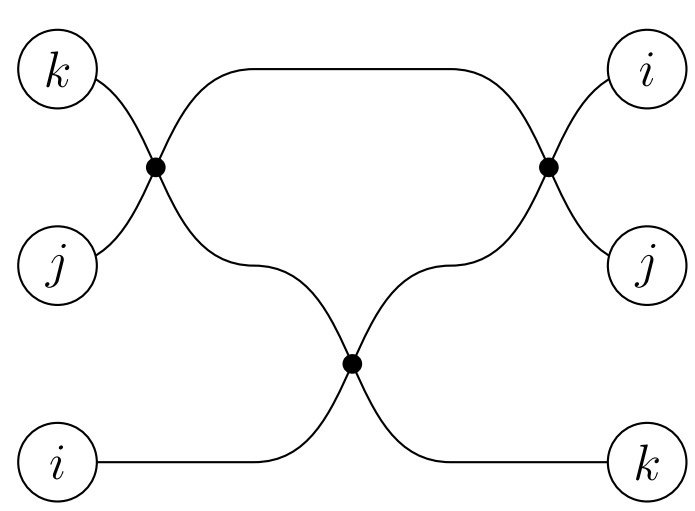

$=$

$=$ $(t_i z_j + z_i)(t_i z_k + z_i)(t_j z_k + z_j)$ 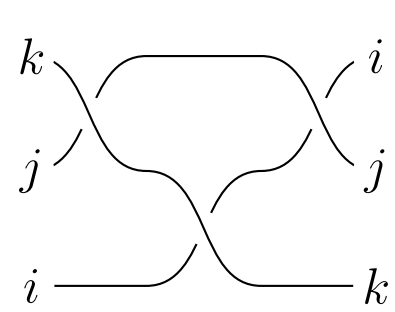

$+$ $(t_i z_j + z_i)(t_i z_k + z_i)$ 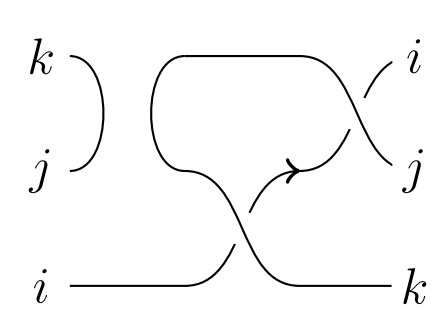

$+ \quad (t_iz_k + z_i)(t_jz_k + z_j)$ 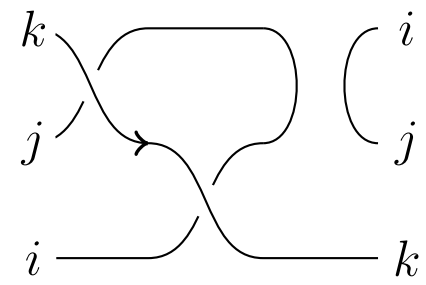 $+ \quad (t_iz_j + z_i)(t_jz_k + z_j)$ 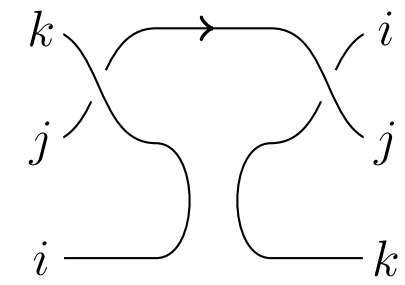

$+ \quad (t_iz_k + z_i)$ 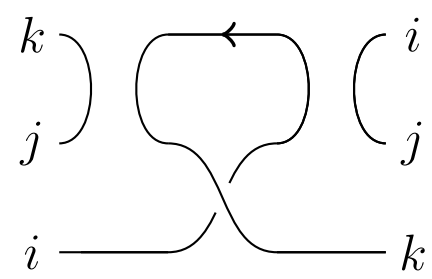 $+ \quad (t_iz_j + z_i)$ 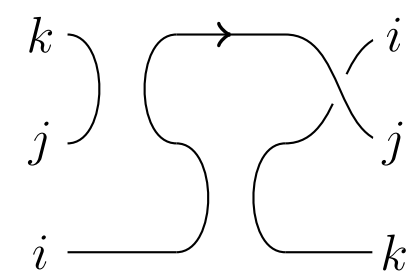

$+ \quad (t_jz_k + z_j)$ 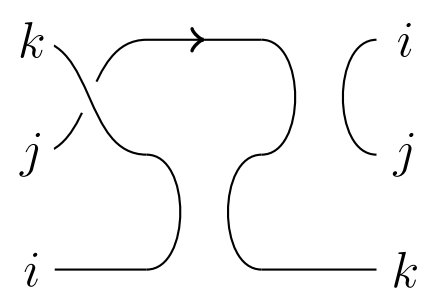 $+$ 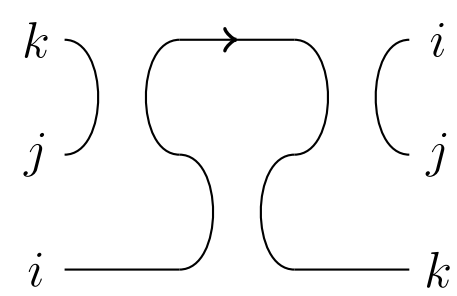

Firstly, notice that the first and fourth term on the left-hand side of the equation cancel with the first and fourth term on the right-hand side of the equation, respectively. Moreover, by Lemma 7.2 we see that the second and third terms on both sides of the equation cancel as well.

Thus, in order to prove the Yang-Baxter equation we have to show that

$(t_iz_k + z_i)$ 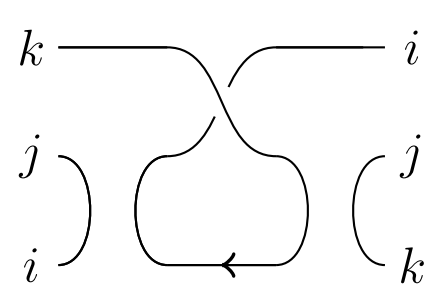 $+ \quad (t_iz_j + z_i)$ 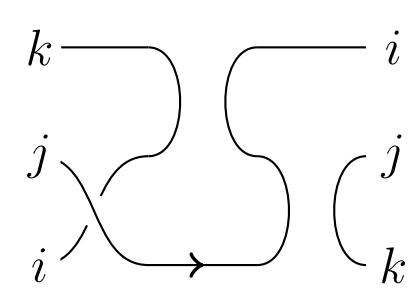

$+ \quad (t_jz_k + z_j)$ 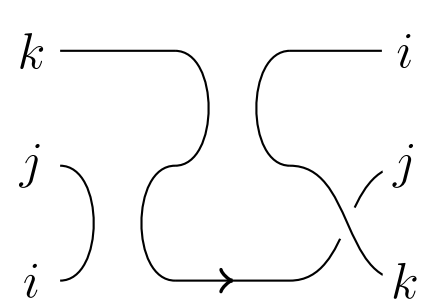 $+$ 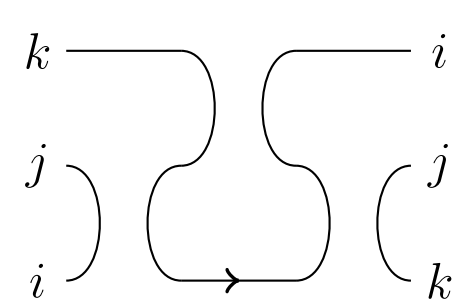 $=$

$= \quad (t_iz_k + z_i)$ 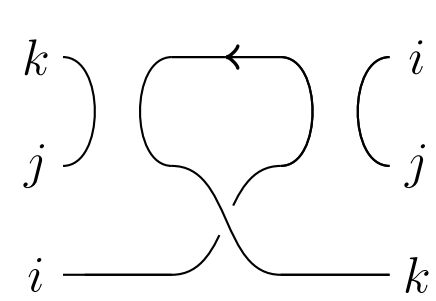 $+ \quad (t_iz_j + z_i)$ 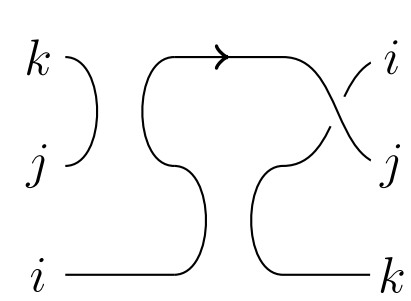

$+ \quad (t_j z_k + z_j)$ 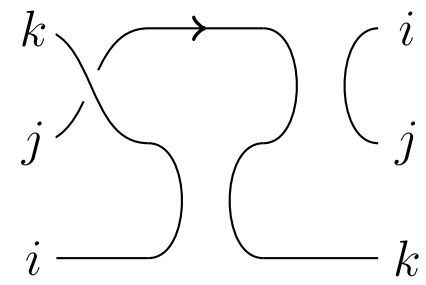 $+$ 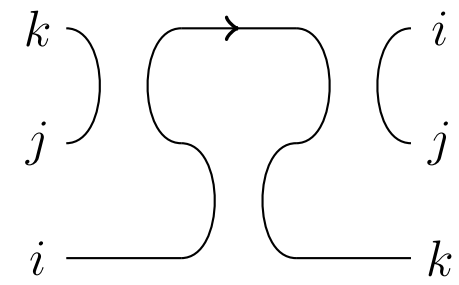

We know the direction of the directed lines going from left to right. Thus, we can move the weights of the caps constituting these directed lines into the respective coefficients for all terms. This allows us to straighten the directed lines and assume they contribute a factor of 1.

$(t_i z_k + z_i)(t_k + 1)z_j$ 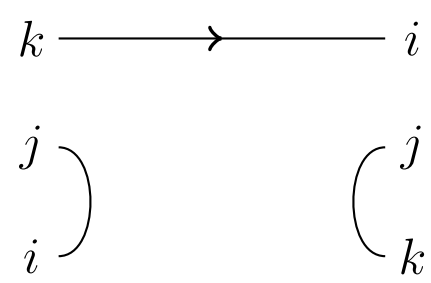 $+ \; (-1)(t_i z_j + z_i)(t_j + 1)z_k$ 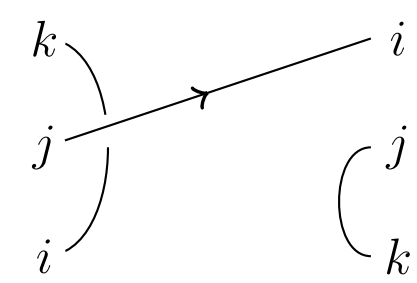

$+ \quad (-1)(t_j z_k + z_j)(t_k + 1)z_i$ 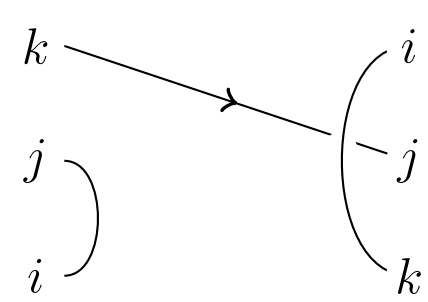 $+ \quad (t_k + 1)z_i(t_j + 1)z_k$ 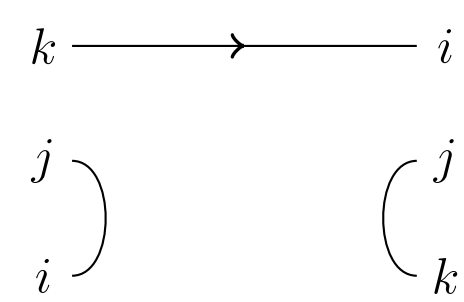 $=$

$= \quad (t_i z_k + z_i)(t_i + 1)z_j$ 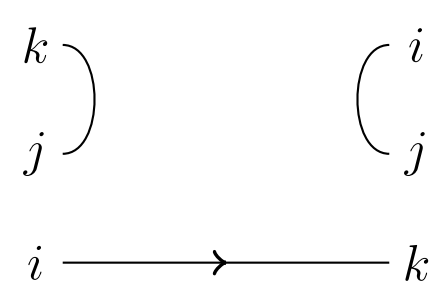 $+ \; (-1)(t_i z_j + z_i)(t_i + 1)z_k$ 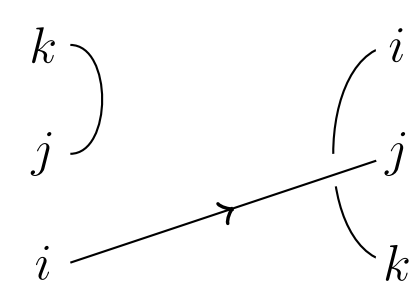

$+ \quad (-1)(t_j z_k + z_j)(t_j + 1)z_i$ 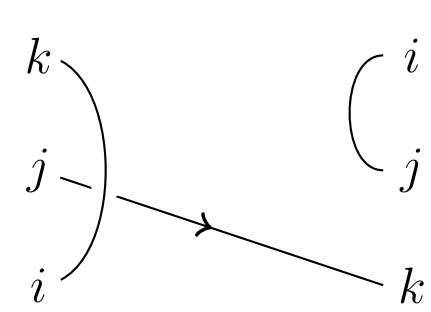 $+ \quad (t_i + 1)z_k(t_j + 1)z_i$ 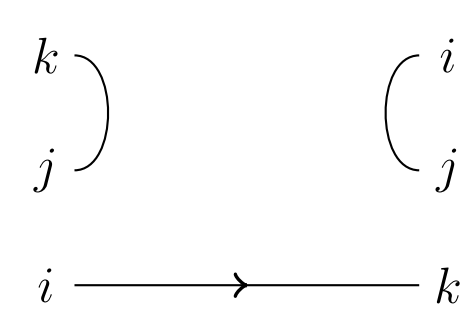

After regrouping terms we get

$[(t_i z_k + z_i)z_j + (t_j + 1)z_i z_k]$ $\Big[$ $(t_k + 1)$ 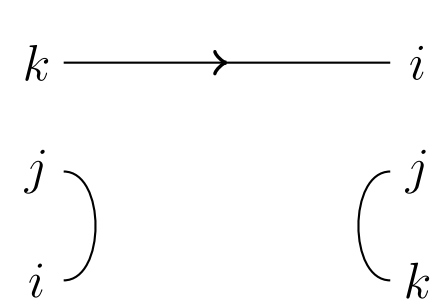 $-$ $(t_i + 1)$ 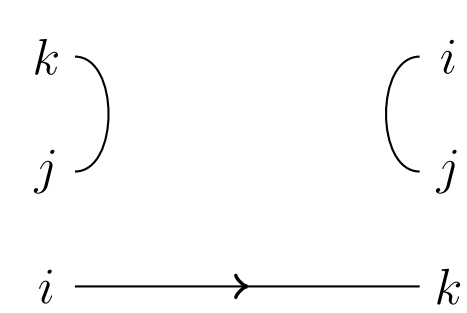 $\Big]$ $=$

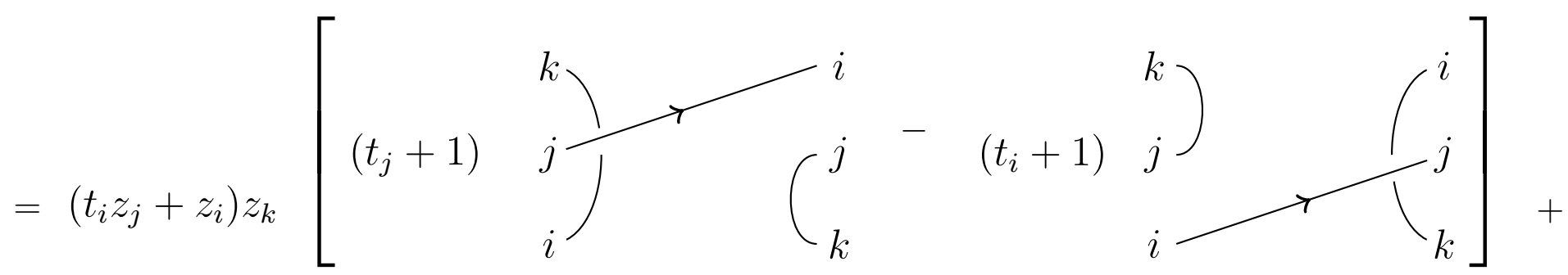

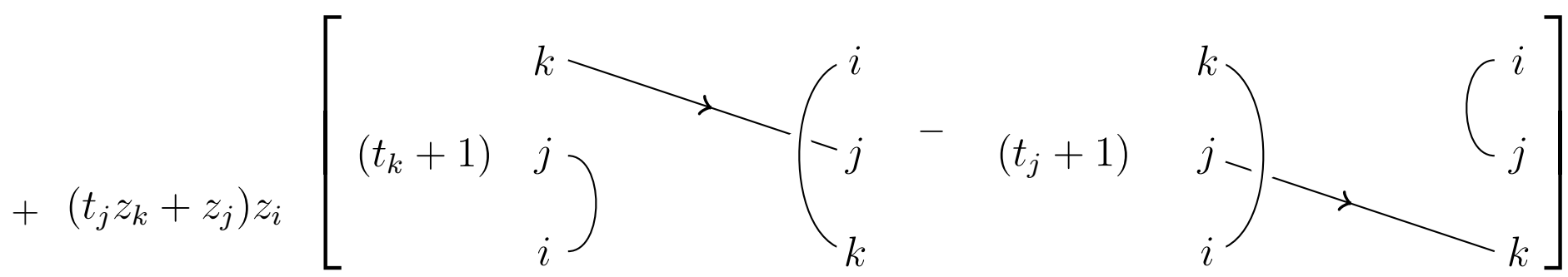

By Corollary 7.5 the three terms in the respective big square brackets are all equal. Hence, the Yang-Baxter equation reduces to checking the identity

$$(t_i z_k + z_i) z_j + (t_j + 1) z_i z_k = (t_i z_j + z_i) z_k + (t_j z_k + z_j) z_i \tag{23}$$

which is true. This concludes the proof of the first case of Theorem 5.1.

**Second case:** This case is almost the same as the first case. So we will skip some steps. Nevertheless, we provide all important parts. Firstly, define the following cap and diagonal line weights:

| | | | |
|---|---|---|---|
| 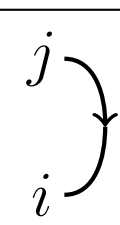 $= z_j$ | 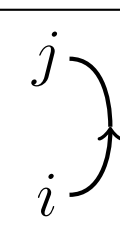 $= -z_i$ | 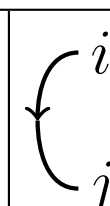 $= (t_i + 1)$ | 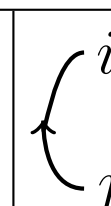 $= (t_j + 1)$ |
| 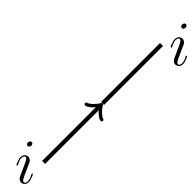 $= -1$ | 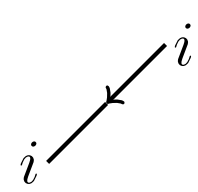 $= 1$ | 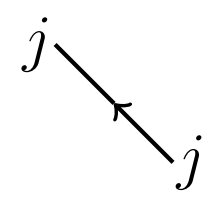 $= 1$ | 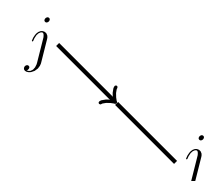 $= -1$ |

Notice that the cap weights are coming from our solutions to the fish equation in Theorems 6.1 and 6.2.

As before the cap and line weights have a few interesting properties.

**Lemma 7.6.** *The Boltzmann weights in Table 8, excluding $a_2$, satisfy the following skein relation*

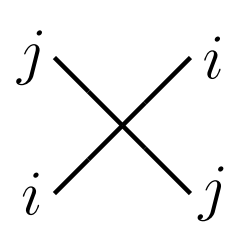 $=$ 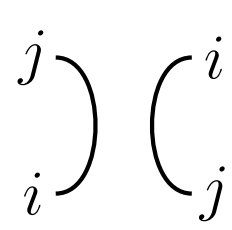 $+ (t_j z_i + z_j)$ 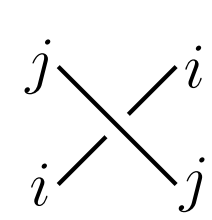 (24)

*where the second term in the skein relation represents the crossing of two straight diagonal lines where the upper line crosses above the lower line.*

*Proof.* The lemma follows by direct computation. □

Next we have the Reidemeister Type 2 move:

**Lemma 7.7.** *The following identity holds*

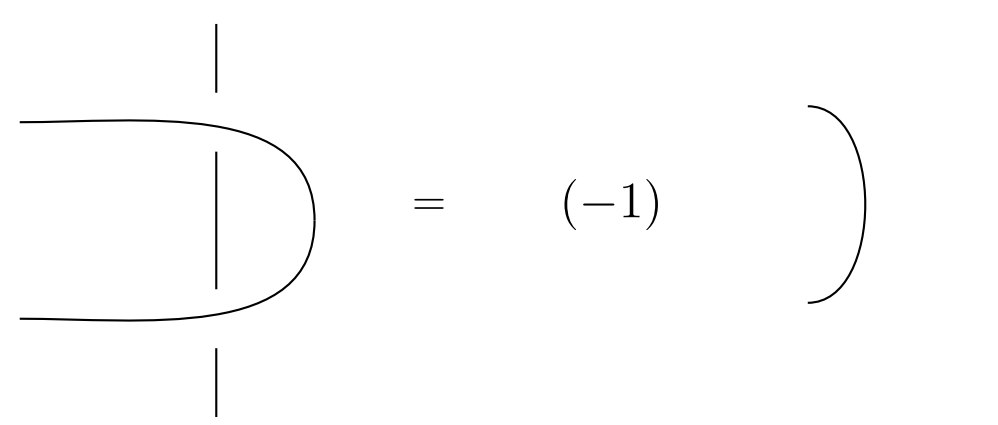

*where the straight line on the right-hand side contributes a factor of* 1.

*Proof.* The −1 factor appears on the right-hand side because the big cap on the left will have one crossing that contributes a −1 factor and one crossing that contributes a 1 factor. On the other hand, the two crossing of the straight line either contribute two −1 factors or two 1 factors. □

Finally, we have the following intriguing identities which come up later in our proof.

**Lemma 7.8.** *The following identities hold*

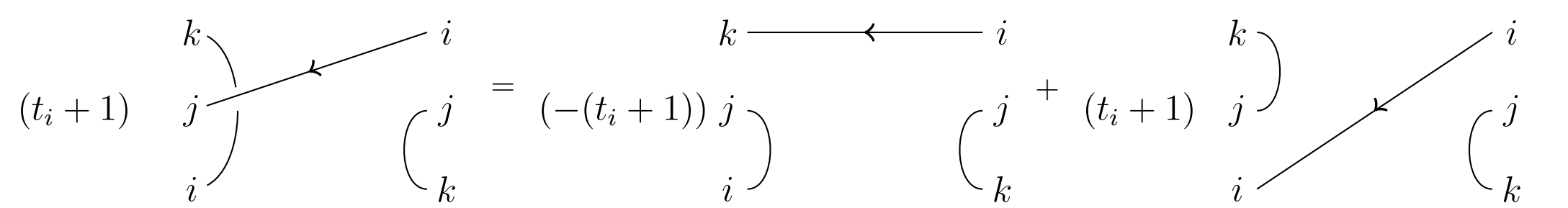

(25)

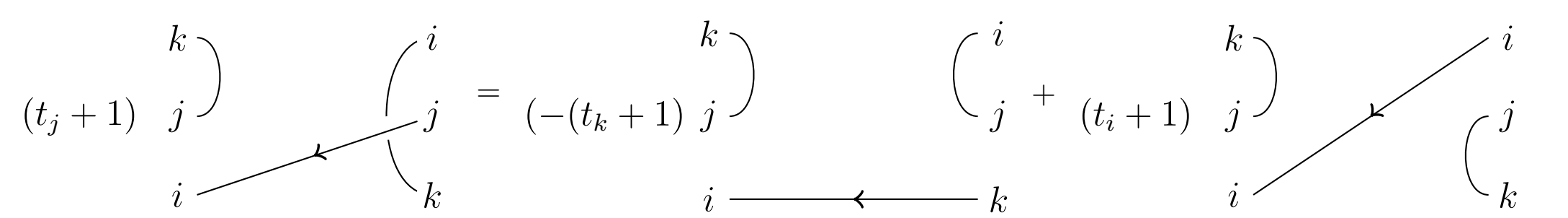

(26)

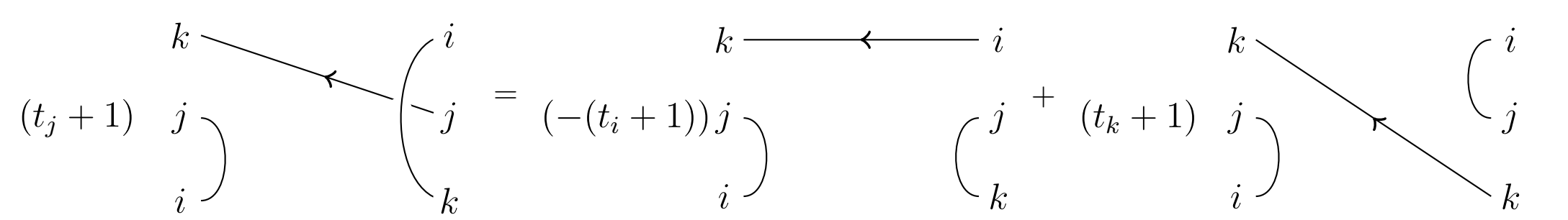

(27)

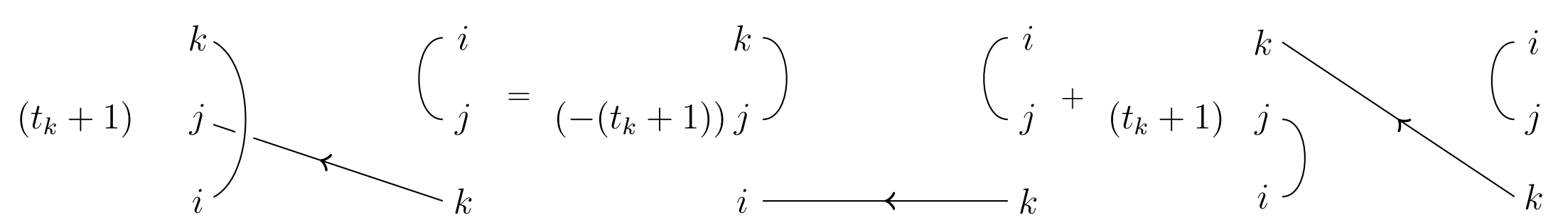

(28)

*where the big caps connecting the rows $i$ and $k$ should be considered as ordinary caps.*

*Proof.* Notice that identities (25) and (28) are equivalent. The same is true for identities (26) and (27). Now identities (25) and (26) are straightforward to check by direct computation. However, when checking the identities, one has to be careful and remember that some crossings can contribute a −1 factor. □

**Remark 7.9.** *The four identities in Lemma 7.8 resemble skein relations. Indeed, notice that if we substitute the crossing on the left-hand side with two vertical caps or two horizontal caps, respectively, we get the shapes on the right-hand side.*

A direct consequence of Lemma 7.8 is the next corollary:

**Corollary 7.10.** *The following triple equality holds.*

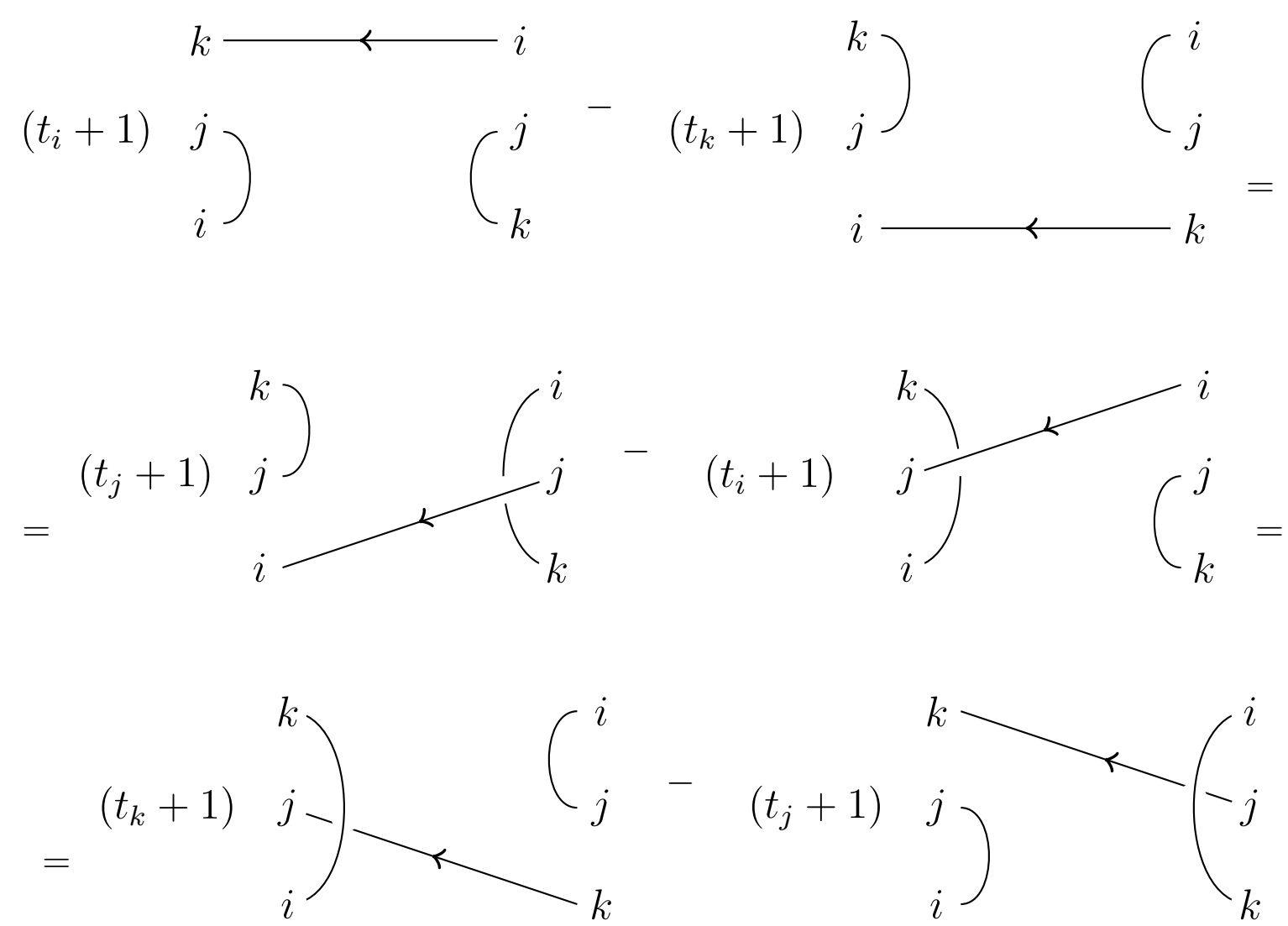

(29)

Now we are ready to prove the Yang-Baxter equation when there are at least four external edges pointing to the left.

By using the skein relation (24) we can expand both sides of the equation into eight terms. We skip the expansion as it is done in the same way as in the first case. Once again, the first four terms on both sides of the equation cancel by using Lemma 7.7. Afterwards, we move the weights of the caps constituting the directed lines from right to left into the coefficients of the respective terms. After regrouping terms we get:

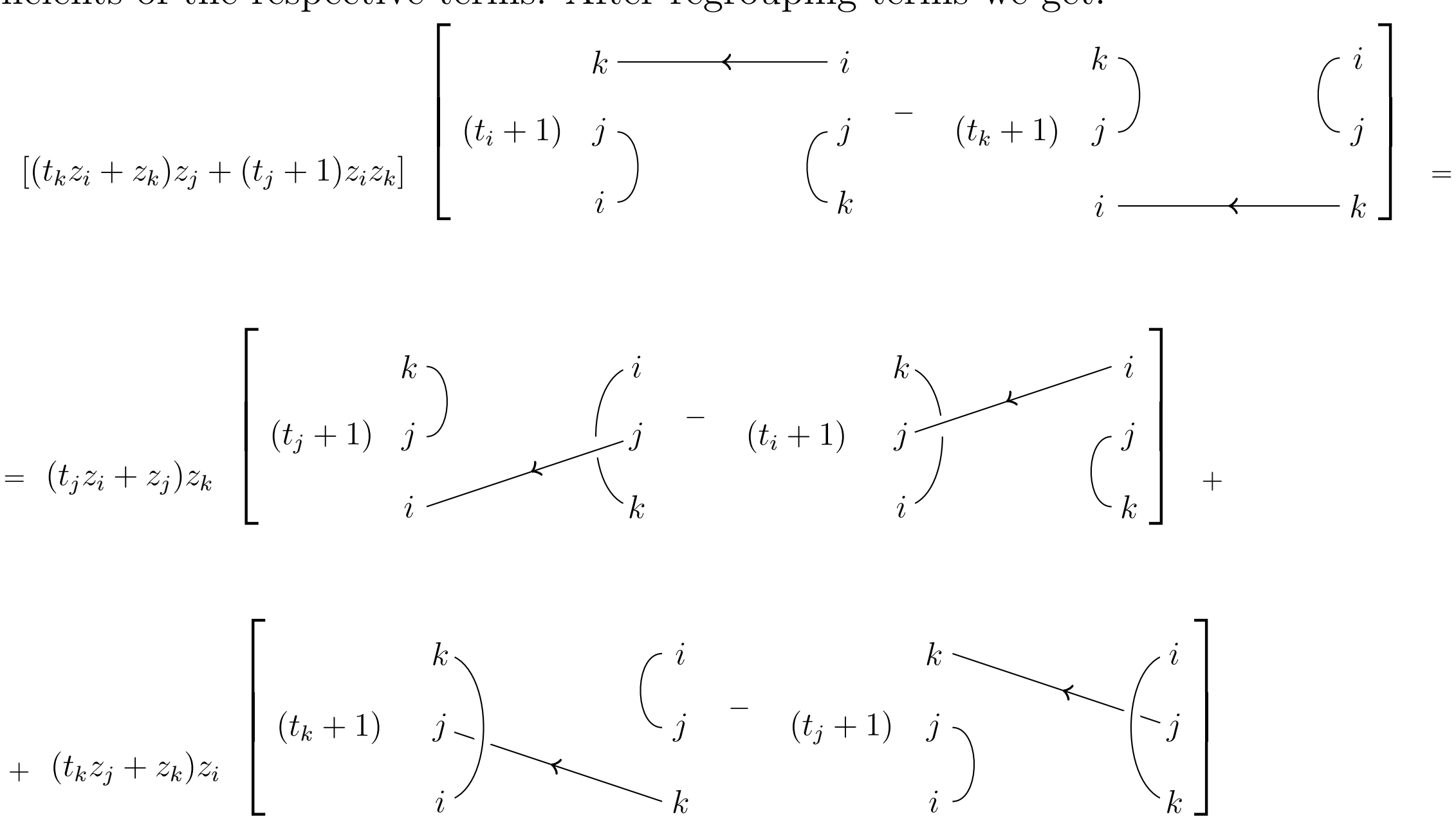

By Corollary 7.10 the three terms in the respective big square brackets are all equal. Hence, the Yang-Baxter equation reduces to checking the identity

$$(t_k z_i + z_k)z_j + (t_j + 1)z_i z_k = (t_j z_i + z_j)z_k + (t_k z_j + z_k)z_i \tag{30}$$

which, remarkably, is exactly identity (23) by switching $i$ and $k$. This concludes the proof of the second case of Theorem 5.1.

Hence, our proof of Theorem 5.1 is complete. □

# Acknowledgements

I am deeply grateful to Prof. Daniel Bump whose constant encouragement and advise throughout all steps of the research project allowed me to keep pushing forward even when I was stuck. This project would not have been possible without him. I would also like to thank Slava Naprienko who gave valuable advice on an earlier draft of this project.

# A Fish equation

The right fish equation is represented by the two equations

$$c_1 \cdot \beta(j,i) + b_2 \cdot \alpha(j,i) = c(i,j) \cdot \beta(i,j) \tag{31}$$

$$c_2 \cdot \alpha(j,i) + b_1 \cdot \beta(j,i) = c(i,j) \cdot \alpha(i,j) \tag{32}$$

These equations can be presented in a matrix form

$$\begin{pmatrix} b_2 & c_1 \\ c_2 & b_1 \end{pmatrix} \begin{pmatrix} \alpha(j,i) \\ \beta(j,i) \end{pmatrix} = c(i,j) \begin{pmatrix} \alpha(i,j) \\ \beta(i,j) \end{pmatrix}$$

Thus, the matrix multiplication switches the indices $i$ and $j$ at the expense of a constant. Let

$$M(i,j) = \begin{pmatrix} b_2 & c_1 \\ c_2 & b_1 \end{pmatrix}$$

Thus,

$$M(j,i) = \begin{pmatrix} -b_2 & c_2 \\ c_1 & -b_1 \end{pmatrix}$$

Hence,

$$M(j,i)M(i,j) \begin{pmatrix} \alpha(j,i) \\ \beta(j,i) \end{pmatrix} = c(i,j)c(j,i) \begin{pmatrix} \alpha(j,i) \\ \beta(j,i) \end{pmatrix}$$

As a result, the matrix $M(j,i)M(i,j)$ has our solutions as an eigenvector with eigenvalue $c(i,j)c(j,i)$.

$$\det[M(j,i)M(i,j)] = \det \left[ \begin{pmatrix} -b_2 & c_2 \\ c_1 & -b_1 \end{pmatrix} \begin{pmatrix} b_2 & c_1 \\ c_2 & b_1 \end{pmatrix} \right] = a_1^2 a_2^2$$

where we used the free-fermionic condition $b_1 b_2 - c_1 c_2 = -a_1 a_2$. Thus, $c(i,j)c(j,i)$ divides $a_1^2 a_2^2 = (t_j z_i + z_j)^2 (t_i z_j + z_i)^2$, which makes it probable that $c(i,j)$ equals $a_1 = (t_j z_i + z_j)$ or $a_2 = (t_i z_j + z_i)$.

In order to find the solutions for the functions $\alpha$ and $\beta$ we just plugged in simple polynomials until we got correct guesses. We expected that there will be 4 solutions as if we accept that $t_i = t_j$, then the equation is easier to work with and we found 4 solutions in that case as well.

The arguments for the left fish equation is analogous.

# References

[B] R. Baxter. Exactly solved models in statistical mechanics. Academic Press Inc. [Harcourt Brace Jovanovich Publishers] London, 1982.

[BIK] N. Bogoliubov, A. Izergin, and V. Korepin. Quantum inverse scattering method and correlation functions. Cambridge University Press, Cambridge, England, 1993.

[BBF] B. Brubaker, D. Bump, S. Friedberg. Schur polynomials and the Yang-Baxter equation, Communications in mathematical physics, Volume 308, Issue 2 (2011) 281–301.

[I1] D. Ivanov, Symplectic ice. In Multiple Dirichlet series, L-functions and automorphic forms, Birkhäuser, Boston, MA. (2012), 205-222.

[I2] A. Izergin. Partition function of a six-vertex model in a finite volume. Soviet Phys. Dokl., 32:878–879, 1987.

[K] G. Kuperberg, Another proof of the alternative-sign matrix conjecture, International Mathematics Research Notices, Volume 1996, Issue 3 (1996), 139–150.

[Z] D. Zeilberger, Proof of the alternating sign matrix conjecture, arXiv preprint math/9407211 (1994)

*E-mail address:* chlalov@stanford.edu